\documentclass[english]{article}

\usepackage[french]{babel}

\usepackage[ansinew]{inputenc}
\usepackage[T1]{fontenc}
\usepackage{makeidx}
\usepackage{setspace}
\usepackage{amssymb}
\usepackage{graphicx}
\usepackage{color}
\makeatletter

\newcommand{\nc}{\newcommand}

\nc{\bat}{|}
\nc{\fty}{_{n \rightarrow \infty}}
\nc{\T}{{\mathbb T}}
\nc{\R}{{\mathbb R}}
\nc{\Q}{{\mathbb Q}}
\nc{\C}{{\mathbb C}}
\nc{\N}{{\mathbb N}}
\nc{\Z}{{\mathbb Z}}
\nc{\K}{{\mathbb K}}
\nc{\DD}{{\mathbb D}}
\nc{\F}{{\mathbb F}}
\nc{\hd}{{\mathcal H}(\DD)}
\nc{\ee}{{\mathcal E}}
\nc{\hw}{H^2_{\sigma}(\mathbb D)}
\nc{\ho}{{\mathcal H}_0(\C \setminus \overline {\mathbb{D}})}
\nc{\hc}{{\mathcal {HF}}(\mathbb T)}
\usepackage{babel}
\makeatother
\makeindex
 \newtheorem{thm}{Theorem}[section]
 \newtheorem{cor}[thm]{Corollary}
 \newtheorem{lem}[thm]{Lemma}
 \newtheorem{prop}[thm]{Proposition}
  \newtheorem{defn}[thm]{Definition}

\date{}

\begin{document}

\title{Do some nontrivial closed $z$-invariant subspaces have the division property?}

\author{Jean Esterle}

\maketitle

Abstract: We consider Banach spaces $E$ of functions holomorphic on the open unit disc $\DD$ such that the unilateral shift $S$ and the backward shift $T$ are bounded on $E.$
Assuming that the spectra of $S$ and $T$ are equal to the closed unit disc we discuss the existence of closed $z$-invariant of $N$ of $E$ having the "division property", which means that the function $f_{\lambda}: z \to {f(z)\over z-\lambda}$ belongs to $N$ for every $\lambda \in \DD$ and for every $f \in N$ such that $f(\lambda)=0.$ This question is related to the existence of nontrivial  bi-invariant subspaces of Banach spaces of hyperfunctions on the unit circle $\T.$ 

\smallskip

\section{Introduction}

Let $\hd$ be the space of holomorphic functions on the unit disc $\DD.$  We denote by $S$ the (unilateral) shift on $\hd$ and by $T$ the "backward  shift" on $\hd,$ defined for $f\in \hd, z \in \DD,$ by the formulae $$Sf(z)=zf(z), \  Tf(z)={f(z)-f(0)\over z},$$

so that $TSf=f$ and $STf=f - f(0).1$ for $f \in \hd.$ 

A linear space $M\subset \hd$ is said to have the {\it division property} if the function $z \to {f(z)\over z-\lambda}$ belongs to $M$ for every $\lambda \in \DD$ and for every $f \in M$ such that $f(\lambda)=0,$ and $N$ is said to be {\it $z$-invariant} if $S(N)\subset N.$

The purpose of the paper is to discuss the existence of non-trivial closed $z$-invariant subspaces having the division property for Banach spaces $E\subset \hd$ satisfying the conditions

$$S(E)\subset E, T(E)\subset E, \mbox{and} \  \rho({S}_{|_{E}})=\rho({T}_{|_{E}})=1,$$

where $  \rho({S}_{|_{E}})$ and $\rho({T}_{|_{E}})$ denote the spectral radii of ${S}_{|_{E}}$ and ${T}_{|_{E}}$.

This question, which is an open problem, was already discussed by A. Volberg and the author in \cite{ev1} for weighted Hilbert spaces of sequences of the form $$H^2_{\sigma}(\DD)=\left \{ f \in \hd \ | \ \Vert f\Vert^2_{\sigma}:= \left [ \sum \limits _{n=0}^{+\infty}\left  \vert {f^{(n)}(0)\over n!}\right \vert^2\sigma^2(n)\right ]^{1\over 2}<+\infty \right \},$$

for weights $\sigma : \Z^+\to (0,+\infty)$ such that $H^2_{\sigma}(\DD)$ satisfies the conditions above. The map $j: f \to \left ({f^{(n)}(0)\over n!}\right )_{n\ge 0}$ is an isometric isomorphism from the space $H^2_{\sigma}(\DD)$ onto to the (complex) sequence space 

$$\ell^2_{\sigma}(\Z^+):=\left \{ u=(u_n)_{n\ge 0} \ | \ \sum \limits_{n=0}^{+\infty}\vert u_n\vert^2\sigma^2(n)<+\infty\right \},$$

and A. Volberg and the author showed that the weight $\sigma$ has an extension to $\Z$ such that every closed translation invariant subspace $F$ of the sequence space $\ell^2_{\sigma}(\Z)$ is generated by the space $F^+:=F\cap\ell^2_{\sigma}(\Z^+),$ which is a closed subspace of $\ell^2_{\sigma}(\Z^+)$ such that $j^{-1}(F^+)$ is a closed $z$-invariant subspace of $\hw$ having the division property. Thus the existence of a weighted Hardy space $\hw$ without nontrivial closed $z$-invariant subspace having the division property would provide an example of an operator $U$ on the separable Hilbert space $H$ such that $U$ and $U^{-1}$ have no common nontrivial closed invariant subpace of $H.$

The motivation of the present paper is to have a fresh look at the existence of non-trivial closed $z$-invariant subspaces having the division property in the context of Banach spaces, which seems to be the natural context to study this question, since the Fr\'echet space $\hd$ provides an obvious example of a Fr\'echet space $E\subset \hd$ without any nontrivial zero- free closed $z$-invariant subspace satisfying

$$ S(E)\subset E, \ T(E) \subset E, \mbox{and} \ \ \rho({S}_{|_{E}})=\rho({T}_{|_{E}})=1.$$

A brief "state-of-the-art" on this question, which was discussed in \cite{ev1}, Sec. 5 for weighted $H^2$-spaces of holomorphic functions on the unit disc, is given in Section 4 of the paper.  

We also discuss the link between closed $z$-invariant subspaces of Banach spaces $E\in \mathcal E$ and closed bi-invariant subspaces of a natural class $\mathcal F$ of Banach spaces of hyperfunctions on the unit circle. Denote by $\ho$ the space of holomorphic functions on $\C \setminus \DD$ vanishing at infinity, and denote by $\hc$ the space of hyperfunctions on $\T$, i.e the set of   pairs $(g,h)$, where $g \in \hd$ and $h\in \ho.$ The nuclear Fr\'echet space $\hc$ can be identified in a natural way to the dual space of the space $\mathcal O(\T)$ of germs of holomorphic functions on $\T,$ and for $f=(f^+,f^-)\in \hc$ the Fourier coefficients of $f$ are defined by the formulae

$$f^+(\zeta)=\sum \limits_{n=0}^{+\infty}\widehat f(n)\zeta^n \ \ (\vert \zeta \vert <1), f^-(\zeta)=\sum \limits_{n<0}\widehat f(n)\zeta^n \ \ (\vert \zeta \vert >1).$$

The (bilateral) shift $\bf S$ on $\hc$ is defined by the formula $$\widehat {{\bf S}f}(n)=\widehat f (n-1) \ (n \in \Z, f \in \hc).$$

Denote by $P^+: (g,h)\to g$ (resp. $P^-:(g,h)\to h)$ the projection of $\hc$ onto $\hd$ (resp. $\ho$). The class $\mathcal F$ is defined to be the family of all Banach spaces continuously contained in $\hc$ such that $P^+(F)\subset F$ which are invariant for $\bf S$ and $\bf S^{-1}$ and for which the spectrum of $\bf S_{|_F}$ is contained in $\T.$ If $F \in \mathcal F,$ then $F\cap \hd=P^+(F)\in \mathcal E$, and, similarly, $F\cap \ho=P^-(F)$ belongs to the class $\mathcal E^-$ of Banach spaces contained in $\ho$ defined in a way similar to the definition of $\mathcal E.$ Conversely
$E\oplus E^-\in \mathcal F$ if $E \in \mathcal E$ and $E^- \in \mathcal E^-.$

We give at the end of the paper two kinds of results, given $E \in \mathcal E.$
\begin{enumerate}

\item Find sufficient conditions on a "tail" $E^- \in \mathcal E^-$ which ensure that  for every closed invariant space $M$ of $E$ having the division property we have

$$\left [ \cup_{n \le 0}{\bf S}^{-n}(M)\right ]^-\cap E=M,$$

where the closure is taken in $E\oplus E^-.$

\item Find a "tail" $E^-\in \mathcal E^-$ such that every closed subspace $N$ of $E\oplus E^-$ such that $\mathcal S^{-1}(N) \subset N$ has the form

$$N=\left [\cup_{n\le 0}{\mathcal S}^{n+k}(M)\right ]^-,$$

where $k\ge 0,$ where $M$ is a closed subspace of $E$ having the division property, and where the closure is taken again in $E\oplus E^-.$
\end{enumerate}
For weighted $H^2$-spaces of holomorphic functions on the disc and weighted $L^2$-spaces of hyperfunctions on the circle, these questions were discussed by A.Volberg and the author in \cite{ev1}, \cite{ev2}. The results of the present paper are analogous to those in \cite{ev1}, \cite{ev2}, but somewhat easier to state in the Banach space context: set $L_n(f)=\widehat f(n)$ for $f\in E\oplus E^-, n \in \Z.$ If the sequence $(\Vert L_{-n}\Vert)_{n\ge 1}$ decreases sufficiently quickly as $n\to +\infty$, then property 1) holds for every nontrivial closed subspace $M$ of $E$ having the division property, and these rates of decrease can be evaluated with respect to the sequence $\left ( \Vert T_{|_E}^n\Vert\right )_{n\ge1}$ via the Matsaev-Mogulski estimates \cite{mm} of the growth as $\vert z \vert \to 1^-$ of a function $f={f_1\over f_2}$ where the holomorphic functions $f_1$ and $f_2$ are allowed to have zeroes in $\DD.$ Concerning question 2 we adapt the methods of \cite{ev2}, which are based on the notion of Dynkin transform and on lower estimates as $\vert z \vert \to 1^-$ for asymptotically holomorhic functions on $\DD$ to show that weighted Hilbert spaces $E^-=H^2_{0,\sigma}(\C \setminus \overline \DD)$ satisfy property 2 with respect to $E$ if the sequence $(\sigma(-n))_{n\ge 1}$ grows sufficiently quickly and regularly as $n\to +\infty.$

The author hopes that the present paper will encourage experts in complex analysis and operator theory to consider this intriguing question of existence of nontrivial closed invariant subspaces having the division property which, besides its link with the bi-invariant subspace problem, presents a lot of interest in itself.
\section{The division property}

We begin this section by recalling some standard properties the unilateral shift $S$ and the backward shift $T$ on the Fr\'echet algebra $\hd.$

We have, or $f\in \hd, p\ge 0, z \in \DD,$ with the convention $\widehat f(n):={f^{(n)}(0)\over n!}$ for $n \ge 0,$ $$T^pf(z)=\sum \limits_{n=p}^{+\infty}\widehat f(n)z^{n-p}, S^{p}T^pf(z)=\sum \limits_{n=p}^{+\infty}\widehat f(n)z^{n},$$ $$(f-S^{p+1}T^{p+1}f )(z)=\sum \limits_{n=0}^p\widehat f(n)z^n, \mbox{and} \ \widehat f(n).1=T^pf-ST^{p+1}f.$$

Let $f \in \hd,$ let $s \in (0,1),$ let $r \in (s, 1), $ and set $M_{r}=\max_{\vert z \vert ={r}}\vert f(z)\vert.$ It follows from Cauchy's inequalites that we have, for $n \ge 0,$

$$\left \vert \widehat f(n)\right \vert \le r^{-n}M_{r},$$
and so, for $\vert z \vert \le s,$

$$\vert (T^pf)(z)\vert \le r^{-p}M_{r}\sum \limits_{n=0}^{+\infty}r^{-n}\vert z\vert^n={r^{-p}M_r\over 1- r^{-1}s}.$$

So the series $\sum \limits_{p=0}^{\infty}\lambda^pT^pf$ converges uniformly on compact subsets of $\DD,$ and $I-\lambda T$ is invertible for $\lambda \in \DD.$

We have
\begin{equation} T(I-\lambda T)^{-1}(S-\lambda I)=(I-\lambda T)^{-1}T(S-\lambda I)=(I-\lambda T)^{-1}(I-\lambda T)=I.\end{equation}

\begin{equation}f(\lambda).1=\sum \limits_{p=0}^{+\infty}\widehat f(p)\lambda^p.1=(I-ST)\left (\sum \limits_{p=0}^{+\infty}\lambda ^pT^p f\right ) =(I-ST)(I-\lambda T)^{-1}f.\end{equation}

Define $f_{\lambda}$ for $f \in \hd, \lambda \in \DD$ by the formula

\begin{equation}f_\lambda(z)={f-f(\lambda)\over z-\lambda} \ \ (z\in \DD).\end{equation}

We have

$$(S-\lambda I)f_{\lambda}=f-f(\lambda.1)=(I-\lambda T)(I-\lambda T)^{-1}f-(I-ST)(I-\lambda T)^{-1}f$$ $$=(S-\lambda I)T(I-\lambda T)^{-1}f.$$ Since the map $f\to (S-\lambda I)f$ is one-to-one on $\hd,$ we obtain

\begin{equation}f_{\lambda}=T(I-\lambda T)^{-1}f \ \ (f\in \hd, \lambda \in \DD).\end{equation}

Obviously, $Spec(S)=\DD,$ since $(S-\lambda I)(\hd)\subsetneq \hd$ for $\lambda \in \DD$ and since the function $z \to {1\over z-\lambda}$ belongs to $\hd$ for $\vert \lambda \vert \ge 1.$ Also if we set $u_\lambda(z)={1\over 1-\lambda z}$ for $\lambda \in \overline{\DD}, z \in \DD,$ then $u_{\lambda}\in Ker(T-\lambda I),$ and so $Spec(T)=\overline \DD$ since $I-\lambda T$ is invertible for $\lambda \in \DD.$

For $f \in \hd, $ set $Z(f):=\{ \lambda \in \DD \ | \ f(\lambda)=0\},$ and for $M \subset \hd,$ set $Z(M):=\cap_{f\in M}Z(f).$ We have the following obvious observation.

\begin{prop} Let $M \neq \{0\}$ be a linear subspace of $\hd,$ and let $\lambda \in \DD.$ Then the following conditions imply each other

(i) $\lambda \notin Z(M),$ and $dim\left (M/((S-\lambda I)(M)\cap M\right ))=1.$

(ii) the function  $f_{\lambda}$ belongs to $M$ for every $f \in M$ such that $f(\lambda)=0.$

\end{prop}

Proof: Assume that (i) holds. Since $\lambda \notin Z(M),$ there exists $u \in M$ such that $u(\lambda)=1,$ and since  $dim\left (M/((S-\lambda)(M)\cap M\right ))=1,$ there exists for every $f \in M$ a function $g \in M$ and a complex number $\mu$ such that $g \in (S-\lambda)(M)$ and $f=\mu u +g.$ Since $\lambda \in Z(g),$ we have $\mu =f(\lambda),$ and so there exists $h \in M$ such that $f(z)=f(\lambda)u +(z-\lambda)h(z)$ for $z \in \DD.$ This shows that $f_{\lambda}\in M$ if $f(\lambda)=0,$ and (i) implies (ii).

Now assume that (ii) holds, and let $g \in M\setminus \{0\}.$ It follows from (ii) that there exists $k \ge 0$ and $h \in M$ such that $g=(S-\lambda I)^kh,$ with $h(\lambda) \neq 0.$ Hence $\lambda \notin Z(M),$ and setting $u:=h(\lambda)^{-1}h,$ we obtain again $u \in M$ such that $u(\lambda)=1.$

Let $f\in M.$ Then $\lambda \in Z(f-f(\lambda ) u),$ and $f=f(\lambda)u +(S-\lambda I)(f-f(\lambda) u)_{\lambda},$ and (ii) implies (i). \ \ \  $\square$

\begin{defn} A linear subspace $M$ of $\hd$ is said to be $z$-invariant if $S(M)\subset M.$ 

\end{defn}

\begin{defn}  A linear subspace $M$ of $\hd$ is said to have the division property at $\lambda \in \DD$ if $f_{\lambda}\in M$ for every $f \in M$ such that $f(\lambda)=0,$ and $M$ is said to have the division property if $M$ has the division property at $\lambda$ for every $\lambda \in \DD.$
\end{defn}

\begin{cor} A linear subspace $M\neq \{0\}$ of $\hd$ has the division property  if and only if it satisfies the two following conditions

(i) $Z(M)=\emptyset$

(ii)  $dim\left (M/((S-\lambda I)(M)\cap M\right ))=1$ for every $\lambda \in \DD.$

\end{cor}

We will see in the next section for closed subspaces $M$ of  Banach spaces of holomorphic functions on $\DD$ satisfying some natural conditions that if $Z(M)=\emptyset$ then the fact that $M$ satisfies property (ii) at some $\lambda\in \DD$ implies that $M$ satisfies property (ii) at $\lambda$ for every  $\lambda\in \DD.$ In other terms if $Z(M)=\emptyset$ the fact that $M$ has the division property at $\lambda$ for some $\lambda \in \DD$ implies that $M$ has the division property at $\lambda$ for every $\lambda \in \DD.$

A linear subspace $M$ of a linear space $E\neq \{0\}$ is as usual said to be nontrivial if $M\neq E$ and $M\neq \{0\}.$ We conclude this section with the following algebraic observation.

\begin{prop} Let $E\neq \{0\}$ be a $z$-invariant subspace of $\hd$ such that $T(E)\subset E$ and $Spec(T_{|_E})\subset \overline \DD,$  let $\lambda \in \DD,$ let $M$ be a nontrivial linear subspace of $E$,  and let $\pi: E\to E/M$ be the canonical surjection. 

(i)  $\lambda \notin Z(M)$ if and only if the map $\pi \circ (S_{|_E}-\lambda I_E): E\to E/M$ is onto.

(ii)  $M$ has the division property at $\lambda$  if and only if we have
$$Ker(\pi \circ(S_{|_E}-\lambda I_E))\subset M.$$ In this case $\pi\circ (S_{|_E}-\lambda I_E)(E)=E/M$ and if we set $U_{\lambda}(\pi (S_{|_E}f-\lambda f))= \pi(f) $ for $f \in M,$ then $U_{\lambda}: E/M \to E/M$  is well-defined and onto.

(iii) If, further, $M$ is $z$-invariant, then $M$ has the division property at $\lambda$  if and only if $S_M-\lambda I_{E/M}$ is one-to-one, where  $S_M:E/M \to E/M$ is defined by the formula $S_M\circ \pi =\pi \circ S_{|_E},$ and in this situation $S_M-\lambda I_{E/M}$ is a bijection and $U_{\lambda}=(S_M-\lambda I_{E/M})^{-1}.$

(iv) Assume that $M$ has the division property at $\lambda$,  and let $\mu \in \DD.$ Then $M$ has the division property at $\mu$  if and only if 
the map $I_{E/M} -(\mu -\lambda)U_{\lambda}$ is one-to-one. In this case $I_{E/M} -(\mu -\lambda)U_{\lambda}$ is a bijection, and we have

$$U_{\mu}=U_\lambda \circ(I_{E/M}-(\mu -\lambda)U(\lambda))^{-1}.$$

\end{prop}

Proof: Since $Spec(T_{|_E})\subset \overline \DD,$ it follows from (4) that $f_{\lambda}$ in $E$ for every $f \in E,$ and so $Z(E)=\emptyset.$

(i) Assume that $\lambda \notin Z(M),$ let $f \in E,$ and let $u \in M$ such that $u(\lambda)=1.$ Then $\pi(f)=\pi (f-f(\lambda)u)=(\pi \circ (S_{|_E}-\lambda I_E)((f- f(\lambda)u)_{\lambda})$ and so $\pi \circ (S_{|_E}-\lambda I_E)$ is onto.

Conversely if $\pi \circ (S_{|_E}-\lambda I_E)$ is onto, let $f \in E$ such that $f(\lambda) =1.$ Then there exists $g \in E$ and $h \in M$ such that $f=(S_{|_E}-\lambda I_E)g+h,$ and $h(\lambda)=1,$ so that $\lambda \notin Z(M).$

(ii) Assume that $Ker(\pi \circ(S_{|_E}-\lambda I_E))\subset M,$ and let $f \in M$ such that $f(\lambda)=0.$ Then $0=\pi(f)=\pi((S_{|_E}-\lambda I_E)f_{\lambda}),$ $f_\lambda \in Ker(\pi \circ(S_{|_E}-\lambda I_E)),$ and so $f_\lambda \in M,$ which shows that $M$ has the division property at $\lambda.$

Conversely assume that $M$ has the division property at $\lambda.$ 

Let $g \in Ker(\pi \circ(S_{|_E}-\lambda I_E)),$ and set $f= (S_{|_E}-\lambda I_E)g.$ Then $f\in M, f(\lambda)=0,$ and so $g =f_{\lambda}\in M.$

In this situation $\lambda \notin Z(M)$ and  it follows from (i) that $\pi\circ (S_{|_E}-\lambda I_E)(E)=E/M.$ For $f \in E,$ set

$$U_{\lambda}(\pi \circ(S_{|_E}-\lambda I_{E})(f))=\pi(f).$$

Then $U_{\lambda}:E/M \to E/M$ is well-defined, since $Ker(\pi \circ(S_{|_E}-\lambda I_E))\subset M,$ and it follows from the definition of $U_\lambda$ that $U_\lambda$ is onto.  

(iii) Now assume that $M$ is $z$-invariant. Then $$\pi\circ (S_{_{|_E}}-\lambda I_E) =(S_M-\lambda I_{E/M})\circ \pi,$$ and it follows from (ii) that $M$ has the division property at $\lambda$ if and only if $S_M-\lambda I_{E/M}$ is one-to-one. In this case it follows from (i) that $S_M-\lambda I_{E/M}$ is also onto, so $S_M-\lambda I_{E/M}:E/M \to E/M$ is a bijection, and we have $U_\lambda = (S_M-\lambda I_{E/M})^{-1}.$ 

(iv) Now assume that $M$ has the division property at $\lambda,$ and let $\mu \in \DD.$ We have

$$(I_{E/M}-(\mu -\lambda)U_{\lambda})\circ \pi \circ (S_{|_E}-\lambda I_E)=\pi\circ (S_{|_E}-\lambda I_E) -(\mu - \lambda)\pi= \pi \circ(S_{|_E}-\mu I_E).$$

If $I_{E/M}-(\mu -\lambda)U_{\lambda}$ is one-to-one, then 

$$Ker( \pi \circ(S_{|_E}-\mu I_E))\subset Ker( \pi \circ (S_{|_E}-\lambda I_E)) \subset M,$$ and so it follows from (ii) that $M$ has the division property at $\mu.$

Conversely assume that $M$ has the division property at $\mu \in \DD,$  and let $u \in Ker(I_{E/M}-(\mu -\lambda)U_{\lambda}).$ There exists $f \in E$ such that $u=\pi(S_{|_E}f-\lambda f),$
and $f \in Ker(\pi\circ (S_{|_E}-\mu I_E)\subset M.$ So $S_{|_E}f \in \mu f +M \in M, u =0,$ and $I_{E}-(\mu - \lambda)U_{\lambda}$ is one-to-one. Also $\mu \notin Z(M),$ $\pi \circ(S_{|_E}-\mu I_E)$ is onto, $I_{E/M}-(\mu -\lambda)U_{\lambda}$ is onto and $I_{E/M}-(\mu -\lambda)U(\lambda)$ is in fact bijective. We have

$$U_{\mu}\circ (I_{E/M}-(\mu -\lambda)U_{\lambda})\circ \pi \circ (S_{|_E}-\lambda I_E)=U_{\mu}\circ (S_{|_E}-\mu I_E)=\pi=U_{\lambda}\circ \pi \circ (S_{|_E}-\lambda I_E).$$

Since $\pi \circ (S_{|_E}-\lambda I_E): E\to E/M$ is onto, this gives $$U_{\mu}\circ (I_{E/M}-(\mu -\lambda)U_{\lambda})=U_{\lambda}, \ U_{\mu}=U_\lambda \circ(I_{E/M}-(\mu -\lambda)U_{\lambda})^{-1}.$$

$\square$

Notice that if $M$ is $z$-invariant and has the division property at $\lambda$ and $\mu,$  then $S_M-\lambda I_{E/M}$ and $S_M-\mu I_{E/M}$ are bijective, $U_\lambda=(S_M-\lambda I_{E/M})^{-1}, U_\mu=(S_M-\mu I_{E/M})^{-1},$ and we obtain directly

$$ (S_M-\lambda I_{E/M})^{-1}(I_{E/ M}-(\mu-\lambda)(S_M-\lambda I_{E/M})^{-1})^{-1}=(S_M-\lambda I_{E/M}-(\mu-\lambda)I_{E/M})^{-1}$$ $$=(S_M-\mu I_{E/M})^{-1},$$

which gives the formula of (iv) when $M$ is $z$-invariant.

Notice also that the argument used to prove (iv) circumvents a small gap in the proof of lemma 2.2 of \cite{ev1} where it is implicitely assumed that $U_{\lambda}$ is one-to-one, which is indeed not true if $M$ is not $z$-invariant. This argument can be used to fix the proof of lemma 2.2 of \cite{ev1}.

\begin{cor} Let $E\neq \{0\}$ be a $z$-invariant subspace of $\hd$ invariant for $T$ such that $Spec(T_{|_E})\subset \overline \DD,$ let $M$ be a nontrivial linear subspace of $E$ having the division property,  and let $\pi: E\to E/M$ be the canonical surjection. 

Then for every $\lambda \in \DD$ there exists an onto linear map $U_{\lambda}:E/M \to E/M$ such that $U_{\lambda}\circ \pi \circ (S_{|_F}-I_{E/M})=\pi,$ the map $I_{E/M}-(\mu -\lambda)U_{\lambda}:E/M \to E/M$ is a bijection for every $\mu \in \DD,$ and we have

\begin{equation}U_{\mu}=U_{\lambda}(I_{E/M} -(\mu - \lambda)U_{\lambda})^{-1}.\end{equation}

\end{cor}

\section{A class of Banach spaces of holomorphic functions on the unit disc}

The spectrum and the spectral radius of a bounded linear operator $R$ on a topological linear space $E$ will be denoted $Spec(E)$ and $\rho(E) .$  We will denote by $\mathcal H(\overline \DD)$ the space  of all  holomorphic functions on $ \DD$ which admit an holomorphic extension to the disc $\DD_r:=\{ z \in \C \ | \ \vert z \vert >r\}$ for some $r>0.$

We will consider the following natural class of Banach spaces of holomorphic functions on the open unit disc.

\begin{defn} We denote by $\mathcal E$ the class of Banach spaces $E\subset \hd$ satisfying the following conditions

 $$ (i) \ S(E)\subset E, T(E)\subset E, S_{|_E}: E \to E \ \mbox{and} \ T_{|_E}: E \to E \ \mbox{are continuous},$$  $$ (ii) \ \ \ \ \ \ \ \ \ \ \ \ \ \ \ \ \ \ \ \ \ \  \max(\rho(S_{|_E}),\rho(T_{|_E})\le 1. \ \ \ \ \ \ \ \ \ \ \  \ \ \ \ \ \ \ \ \ \ \ \ \ \ \ \ \ \ \ \ \ \  \ \ \ \ \  \ $$

\end{defn}

Notice that it follows from the closed graph theorem that if $E \subset \hd$ is a Banach space, if $S(E)\subset E, T(E)\subset E,$ and if the injection $j: E \to \hd$ is bounded, then the maps $S_{|_E}: E \to E$ and $ T_{|_E}: E \to E $ are continuous.

\begin{prop} Let $E \in \mathcal E.$ Then $f_{\lambda}\in E$ for $f \in E, \lambda \in \DD,$ the map $f\to f_{\lambda}$ is bounded on $E$ for $\lambda \in \DD,$ $Z(E)=\emptyset,$ the injection $j: E \to \hd$ is bounded, $\mathcal H(\overline{\DD})\subset E$ and $Spec(S_{|_E})=Spec(T_{|_E})=\overline \DD.$

\end{prop}

Proof: Denote by $I_E:f \to f$ the identity map on $E.$ Since $I_E-\lambda T_{|_E}$ is invertible for $\lambda \in \DD,$ it follows from (4) that $f_{\lambda} \in E$ for $f \in E, \lambda \in \DD,$ and that the map $f \to f_{\lambda}$ is bounded on $E$ for $\lambda \in \DD,$ $E$ has the division property, and $Z(E)=\emptyset.$  Since $f-f(0).1= ST(f)\in E$ for $f \in E,$ $E$ contains the constant functions. 

It follows from (2) that the evaluation functional $\delta_{\lambda}: f \to f(\lambda)$ is continuous for $\lambda \in \DD$ and that the injection $j:E \to \hd$ is bounded. Since $E$ has the division property, and since $S(E)\subset E,$ we have $(S-\lambda I_E)(E) =Ker(\delta_{\lambda})$ for $\lambda \in \DD,$ and $Spec(S_{|_E})=\overline{\DD}.$

Let $f \in \mathcal H(\overline{\DD}).$ Since $\rho(S_{|_E})= 1,$ the series $\sum \limits_{n=0}^{+\infty}\widehat f(p)S_{|_E}^p.1$ converges in $E,$ and  $f =\sum \limits_{n=0}^{+\infty}\widehat f(p)S_{|_E}^p.1$ since $\delta_{\lambda}$ is continuous on $E$ for $\lambda \in \DD.$ This shows that $\mathcal H(\overline{\DD})\subset E$. Set $u_\lambda (z)={1\over 1-\lambda z}.$ Then $u_\lambda \in \mathcal H(\overline{\DD})$ and $Tu_{\lambda}=\lambda u_{\lambda}$ for $\lambda \in \DD.$ Hence $\DD \subset Spec(T_{|_E})$ and $Spec(T_{|_E})=\overline \DD.$   $\square$

Notice that it follows from (4) that the map $\lambda \to f_{\lambda}$ is holomorphic on $\DD$ for $f \in E,$ and that it follows from  (2) that the map $\lambda \to \delta_{\lambda}$ is holomorphic from $\DD$ into the dual space $E^*.$ 

Clearly, the closure of a $z$-invariant subspace of a Banach space $E \in \mathcal E$ is also $z$-invariant. We also have the following easy result.

\begin{prop} Let $M$ be a subspace of a Banach space $E\in \mathcal E.$ If $M$ has the division property, then $\overline M$  also has the division property.

\end{prop}

Proof: Assume that  $M$ has the division property, let $f \in \overline M,$ let $\lambda \in \DD$ such that $f(\lambda)=0,$ let $u\in M$ such that $u(\lambda)=1,$ and let $(f_n)_{n\ge 1}$ be a sequence of elements of $M$ such that $\lim_{n\to +\infty}\Vert f-f_n\Vert=0.$ Then $\lim_{ n \to +\infty}f_n(\lambda)=0,$ and so $\lim_{n\to +\infty}\Vert f-f_n-f_n(\lambda)u\Vert=0.$
Since $\lambda \in Z(f_n-f_n(\lambda)u),$ $(f_n-f_n(\lambda)u)_{\lambda}\in M,$ and $f_{\lambda}=\lim_{n\to +\infty}(f_n-f_n(\lambda)u)_{\lambda}\in \overline M.$   $\square$

It follows from Proposition 2.5 that $\lambda \in \DD \setminus Z(M)$ if a linear space $M\subset \hd$ has the division property at $\lambda \in \DD.$ For closed subspaces $M$ of a Banach space $E\subset \mathcal E,$ a standard conectedness argument gives the following result.

\begin{prop} Let $E \in \mathcal E$ be a Banach space, let $M\neq \{0\}$ be a closed subspace of $\mathcal E$ and let $\Omega(M)$ be the set of all $\lambda \in \DD$ such that $M$ has the division property at $\lambda.$ Then either $\Omega(M)=\emptyset,$ or $\Omega(M)=\DD \setminus Z(M).$

\end{prop}

Proof: Assume that $\lambda_0\in \Omega(M).$  It follows from Proposition 2.5 (iv) that $\lambda \in M$ if and only if $I_{E/M}- (\lambda -\lambda_0)U_{\lambda_0}:E/M \to E/M$ is bijective, and so $\Omega(M)$ is open, since the group of invertible elements of a Banach algebra is open.

Now let $\lambda \in \overline{\Omega(M)}\cap (\DD\setminus Z(M)),$ let $(\lambda_n)_{n\ge 1}$ be a sequence of elements of $\Omega(M)$ such that $\lim_{n\to +\infty}\vert \lambda -\lambda_n\vert=0.$
and let $f \in M$ such that $f(\lambda)=0.$ 

Let $h \in M$ such that $h(\lambda)=1.$ We may assume that $h(\lambda_n)\neq 0$ for $n\ge 1.$ Set $g_n=f_n-f_n(\lambda_n)h(\lambda_n)^{-1}h.$ Then $g_n\in M,$ $g_n(\lambda_n)=0,$ $(g_n)_{\lambda_n}\in M,$ $\lim_{n\to +\infty}\Vert f-g_n\Vert=0,$ and it follows from (4) that $\lim_{n\to +\infty}\Vert f_{\lambda}-(g_n)_{\lambda_n}\Vert =0.$ So $f_\lambda \in M,$ and $\Omega(M)$ is a closed subset of $\DD \setminus Z(M).$ Since $\DD\setminus Z(M)$ is connected, this shows that $\Omega(M)=\DD \setminus Z(M).$ 

 $\square$
 
 For $U \subset \hd$ we will use the notation $zM:=S(M).$ We then deduce from Corollary 2.4 and Proposition 3.4 the following result.

 \begin{cor} Let $E \in \mathcal E$ be a Banach space, and let $M\neq \{0\}$ be a closed subspace of $E.$ Then $M$ has the division property if an only if the two following conditions are satisfied.
 
 \smallskip
 
 (i) $Z(M)=\emptyset.$
 
 \smallskip
 
 (ii) $dim(M/(M\cap zM))=1.$
 
 \end{cor}

 \section{Examples of Banach spaces of holomorphic functions with nontrivial $z$-invariant subspaces having the division property}

 Let $E \in \mathcal E.$ A function $f \in E$ is said to be cyclic if $span\{S^n f : n\ge 1\}$ is dense in $E.$ Since the space $\C[z]$ of polynomial functions on $\DD$ has the division property, it follows from Proposition 3.3 that $\overline{span\{S^n f : n\ge 1\}}$ has the division property if $Z(f)=\emptyset,$ so if $E$ contains a non-cyclic function $f$ without zeroes in $\DD$ then $\overline{span\{S^n f : n\ge 1\}}$ is a nontrivial $z$-invariant subspace having the division property. In particular if the polynomial functions are not dense in $E$ then $\overline{\C[z]} $ is a nontrivial closed $z$-invariant subspace of $E$ having the division property. For example the disc algebra $\mathcal A(\DD)$ is a closed $z$-invariant subspace of $H^{\infty}(\DD)$ having the division property. More generally it follows from the following observation that every Banach algebra of holomorphic functions $B \in \mathcal E$ possesses nontrivial $z$-invariant subspaces having the division property. 
 
 \begin{prop} Let $E \in \mathcal E$ be a Banach space, and for $\lambda \in \DD$ let $\delta_{\lambda}: f \to f(\lambda)$ be the evaluation functional of elements of $E$ at $\lambda.$
 
 If $\liminf_{\vert \lambda \vert \to 1^-}\Vert \delta_{\lambda}\Vert<+\infty,$  there exists $\zeta \in \partial \DD$ and $\phi \in E^*$ such that $\phi(1)=1,$ $\phi\circ S_{|_E}=\zeta \phi,$ and  $Ker(\phi)$ \ is a nontrivial $z$-invariant subspace of $E$ having the division property.
 
 \end{prop}
 
 Proof: Assume that $\liminf_{\vert \lambda \vert \to 1^-}\Vert \delta_{\lambda}\Vert<+\infty,$ and let $(\lambda_n)_{n\ge 1}$ be a sequence of elements of $\DD$ such that $\lim_{n\to +\infty}\vert \lambda_n \vert =1$ and $\limsup_{n\to +\infty} \Vert \delta_{\lambda_n}\Vert <+\infty.$ 
 
 Taking a subsequence if necessary we can assume that the sequence $(\lambda_n)_{n\ge 1}$ has a limit $\zeta \in \partial \DD$ as $n\to +\infty.$ Since bounded subset of $E^*$ are relatively compact with respect to the weak$^*$-topology $\sigma(E^*,E),$ the sequence $(\delta_{\lambda_n})_{n\ge 1}$ has a weak$^*$-cluster point $\phi\in E^*,$ and $\langle 1,\phi \rangle=1$ since $\langle 1, \delta_{\lambda_n}\rangle=1$ for $n\ge 1.$ Let $f \in E.$ There exists a subsequence $(\lambda_{n_p})_{p\ge 1}$ of the sequence $(\lambda_n)_{n\ge 1}$ such that
 
 $$\langle f, \phi \rangle =\lim_{p\to +\infty}\langle f, \delta_{n_p}\rangle =f(\lambda_{n_p}), \langle Sf,\phi \rangle =\lim_{p\to +\infty}\langle Sf, \delta_{\lambda_{n_p}}\rangle=\lim_{p\to +\infty}\lambda_{n_p}f(\lambda_{n_p}),$$ and so $\phi\circ S_{|_E}=\zeta\phi.$
 
 Clearly, $Ker(\phi)$ is $z$-invariant. Now let $\lambda \in \DD,$ and let $f \in Ker(\phi)$ such that $f(\lambda)=0.$  Then $$0=\langle f ,\phi \rangle=\langle Sf_{\lambda}-\lambda f_{\lambda}, \phi \rangle=(\zeta -\lambda)
\langle f_{\lambda}, \phi \rangle, $$

and so $f_\lambda \in Ker(\phi),$ which shows that  $Ker(\phi)$ has the division property. $\square$

If $M$ is a $z$-invariant subspace of $\hd,$ define as usual the index of $M$ by the formula

$$ind(M)=dim(M/zM)\in \{1, 2, \dots, +\infty\},$$ and say that $M\subset \hd$ is {\it zero-free} if $Z(M)=\emptyset.$ If follows from Corollary 3.5 that the nontrivial closed $z$-invariant subspaces of a Banach space $E\in \mathcal E$ having the division property are the zero-free closed $z$-invariant subspaces of index 1. Borichev showed in \cite{b} that zero-free closed $z$-invariant subspaces of arbitrary index $k\in \{2, \dots, +\infty\}$ do exist for a large class of Banach spaces of holomorphic functions in the disc, but his construction does not give zero-free closed $z$-invariant subspaces of index 1. Borichev, Hedenmalm and Volberg showed in \cite{bhv} that noncyclic functions $f$ without zeroes in the disc do exist in "large" weighted Bergman spaces, and the existence of noncyclic functions without zeroes for some other classes of spaces of holomorphic functions in the disc follows from works of Atzmon \cite{a} based on the theory of entire functions and from works of Nikolski \cite{n} based on the so-called "abstract Keldysh method", see \cite{ev1}, section 5. We refer more generally  to Section 5 of \cite{ev1} for a detailed discussion of the existence of nontrivial closed $z$-invariant subsbaces having the discussion property for weighted $H^2$-spaces of holomorphic functions on $\DD.$

If $\Vert f \Vert \le \Vert g \Vert$ for every $f , g \in E$ such that $\vert f(\zeta)\vert \le \vert g(\zeta)\vert$ for every $\zeta \in \DD,$ then $\Vert p(S)\Vert \le \max_{\zeta \in \T}\vert p(\zeta)\vert$ for every polynomial $p$, the shift $S$ satisfies the von Neumann inequality, and the Brown-Chevreau-Pearcy theory \cite{ber}, \cite{ch} and its extensions to general Banach spaces \cite{am} apply and produce a very rich lattice of $z$-invariant subspaces of $E$ (we are in the "easy" situation where the spectrum of $S$ equals the whole of $\DD$). The author was not able so far to use these methods to produce zero-free non $z$-cyclic elements of $E$ or more generally closed $z$-invariant subspaces of $E$ in this situation. In the other direction it does not seem that Read's method construction of a counterexample to the invariant subspace problem and its adaptation to Hilbert spaces \cite{gr}, \cite{re} gives any clue to construct a Banach space $E\in \mathcal E$ without nontrivial closed $z$-invariant  subspaces having the division property.

Notice that the Fr\'echet space $\hd$ has no proper closed $z$-invariant subspace having the division property, and more generally, no proper zero-free closed $z$-invariant subspace. To see this, first notice that a closed $z$-invariant subspace $M$ of $\hd$ is a closed ideal of the Fr\'echet algebra $\hd.$ It follows from standard results about Fr\'echet algebras \cite{mi}  that if $M$ is nontrivial then the unital Fr\'echet algebra $\hd/M$ possesses a character, which means that there exists a character $\chi$ on $\hd$ such that $M\subset Ker (\chi).$ Since polynomials are dense in $\hd$ characters of $\hd$ have the form $\chi_{\lambda}:f \to f(\lambda)$ for some $\lambda \in \DD,$ and so $Z(M)\neq \emptyset$ for every nontrivial closed $z$-invariant subspace of $\hd.$

\section{Banach spaces of hyperfunctions on the unit circle}

Denote by $\T=\partial \DD$ the unit circle, denote by $\ho$ the space of holomorphic functions on $\C \setminus \overline \DD$ vanishing at $\infty,$ and denote by $\hc$ the space of hyperfunctions on $\T,$ i.e. the space of all pairs $f=(f^+,f^-)$ where $f^+\in \hd$ and $f^-\in \ho.$ 

Hyperfunctions on the circle form a flabby sheaf \cite{c},  the notion of support of a distribution can be extended to hyperfunctions, the "product" of two hyperfunctions with disjoint support  vanishes in some natural sense \cite{eg} and several variables extensions of the notion of hyperfunction play a basic role in microlocal calculus. We will not use these saddleties here, but we will use some standard properties of $\hc$ considered as a topological convex linear space.

Set, for $f=(f^+,f^-) \in \hc,$ $$p_n (f)=\max(\sup_{\vert \lambda \vert \le 1-1/n}\vert f^+(\lambda)\vert, \sup_{\vert \lambda \vert \ge 1+1/n}\vert f^-(\lambda)).$$

Then $(\hc, (p_n)_{n\ge 1})$ is a Fr\'echet space. Identifying $g$ with $(g,0)$ for $g\in \hd,$ we can consider $\hd$ as a closed subspace of $\hc.$

Since $\ho$ is isomorphic to the space of functions analytic on $\DD$ vanishing at the origin, and since closed subspaces and products of nuclear Fr\'echet spaces are nuclear,  the fact that $\mathcal H(\mathcal U)$ is nuclear for every open subset $\mathcal U$ of $\C,$ see \cite{p}, theorem 6.4.2, implies  that $\hc$ is a nuclear Fr\'echet space.

For $s \in (0,1)$ set $\mathcal U_s:=\{ \lambda \in \C \ : \  s \le \vert \lambda \vert <{s^{-1}}\},$ and denote by $\mathcal O (\T):=\cup_{s\in (0,1)}\mathcal H(\mathcal U_{s})$ the space of germs of analytic functions on $\T,$ equipped with the usual inductive limit topology, for which the bounded sets are the sets which are contained and bounded in $\mathcal H(\mathcal U_{s})$ for some $s \in (0,1).$  One can identify  $\mathcal O (\T)$ with the dual space of $\hc,$ equipped with the topology of uniform convergence on bounded subsets of $\hc,$ by using for $h \in \mathcal H(\mathcal U_s), 0<s<1$ the formula

$$<f, h>={1\over 2i\pi }\int_{r\T}f(\zeta)h(\zeta)d\zeta+{1\over 2 i \pi}\int_{R\T}f(\zeta)h(\zeta)d\zeta \ \ \ (f\in \hc),$$

 where $s < r< 1<R<s^{-1},$ and where the unit circle $\T$ is oriented counterclockwise, see the details in Chapter 1 of \cite{bg}. 

The Fourier coefficients $\widehat f(n)$ and $\widehat h(n)$ for $f \in \hc$ and $h \in \mathcal H(\mathcal U_s)$ are defined by the formulae

$$f^+(\zeta)=\sum \limits _{n=0}^{+\infty}\widehat f(n) \ (\vert \zeta\vert <1), \ f^-(\zeta) =\sum \limits_{n<0}\widehat f(n)\zeta^n \ (\vert \zeta \vert >1),$$ $$ \ h(\zeta)=\sum_{n=-\infty}^{+\infty}\widehat h(n)\zeta^n \ (\zeta \in \mathcal U_s),$$  which gives

$$<f, h>= \sum_{n\in \Z} \widehat f(n)\widehat h(-n-1).$$

It follows from the standard properties of nuclear Fr\'echet spaces \cite{p}, Theorem  4.4.13 that if we set $<h,\tilde f>=<f, h>$ for $f \in \hc, h \in \mathcal O(\T),$ then the application $f \to \tilde f$ is an isomorphism from $\hc$ onto the dual space of $\mathcal O(\T)$, equipped with the topology of uniform convergence on bounded subsets of $\mathcal O(\T).$

We now introduce the (bilateral) shift ${\bf S}$ which is an extension to $\hc$ the (unilateral) shift $S$ defined on $\hd$, by using the formula

$$\widehat{{\bf S}f}(n)=\widehat {\bf S} (n-1), \ \ f=(f^+,f^-) \in \hc, n \in \Z.$$

This gives

$$({\bf S}f)^+(\zeta)=\zeta f^+(\zeta) +\widehat f (-1).1\ \ (\vert \zeta \vert <1), \ ({\bf S}f)^-(\zeta)=\zeta f^-(\zeta) -\widehat f (-1).1 , \ (\vert \zeta \vert >1),$$

$$({\bf S}^{-1}f)^+(\zeta)= \sum \limits_{n=0}^{+\infty}\widehat f(n+1)\zeta^n=\sum \limits_{n=1}^{+\infty}\widehat f(n)\zeta^{n-1} \ (\zeta \in \DD), $$ 
$$({\bf S}^{-1}f)^-(\zeta)= \sum \limits_{n<0}\widehat f(n+1)\zeta^n={f^-(\zeta)+\widehat f(0)\over \zeta}= {f^-(\zeta)+f^+(0)\over \zeta} \ (\vert \zeta \vert >1).$$
In particular ${\bf S}(g,0)=(Sg,0)$ for $g \in \hd,$ and $({\bf S}^{-1}(0,h))^-(\zeta)={h(\zeta)\over \zeta}$ for $h \in \ho, \vert \zeta \vert >1.$

For $h\in {\mathcal H}({\mathcal U}_r), \zeta \in {\mathcal U}_r,$ set ${\bf S}^*h(\zeta)=\zeta h(\zeta).$ We have

$$<h, {\bf S}f>=<{\bf S}^*h,f> \ \ (f\in \hc, h \in {\mathcal O}(\T)).$$

Obviously, ${\bf S}^* -\lambda I_{{\mathcal O}(\T)}: {\mathcal O}(\T)\to {\mathcal O}(\T)$ is invertible if and only if $\lambda \notin \T,$ and so $Spec({\bf S})=\T=Spec({\bf S}^{-1}).$



Denote by $P^+: (f^+,f^-) \to f^+$ and $P^-: (f^+,f^-) \to f^-$ the projection maps. We now introduce a natural class of Banach spaces of hyperfunctions on the circle.

\begin{defn} Let $\mathcal F$ be the class of (nonzero) Banach spaces $F \subset \hc$ satisifying the following properties

\smallskip

(i) $P^+(F)\subset F,$ and $P^+_{|_F}:F\to F$ is continuous.

\smallskip

(ii)  ${\bf S}(F)\subset F,$ ${\bf S}_{|_F}:F\to F$ is continuous, and $Spec({\bf S}_{|_F})\subset \T.$

\end{defn}
Notice that if the  injection $j:F \to \hc$ is continuous, and if $P^+(F)\subset F$ and ${\bf S}(F)\subset F,$ then it follows from the closed graph theorem that $P^+_{|_F}:F\to F$ and ${\bf S}_{|_F}:F\to F$ are continuous. Also if ${\bf S}_{|_F}:F\to F$ is continuous, and if $0\notin Spec({\bf S}_{|_F}),$ then it follows of course from the closed graph theorem that $({\bf S}_{|_F})^{-1}:F\to F$ is continuous. Notice also that an obvious verification shows that the fact that ${\bf S}_{|_F}$ is invertible is equivalent to the fact that ${\bf S}^{-1}(F)\subset F,$ and that in this case ${\bf S}^{-1}_{|_F}=({\bf S}_{{|_F}})^{-1}.$

We now introduce a natural class of Banach spaces contained in $\ho$.

\begin{defn} For $g\in \ho,$ $\vert \zeta \vert >1,$ set

 $$S^- g(\zeta)={\zeta g(\zeta)}-\widehat f(-1), T^-g(\zeta)={g(\zeta)\over \zeta}.$$  
 
 We denote by $\mathcal E^-$ the class of Banach spaces $E^-\subset {\mathcal H}_0(\C \setminus \DD)$ satisfying the following conditions

 $$ (i) \ S^-(E^-)\cup T^-(E^-)\subset E^-, \ \mbox{and the maps} \ S^-_{|_{E^-}}: E^- \to E^- \ \mbox{and} \ T^-_{|_{E^-}}: E^- \to E^- $$ are continuous,  $$ (ii) \ \ \ \ \ \ \ \ \ \ \ \ \ \ \ \ \ \ \ \ \ \  \max(\rho(S^-_{|_{E^-}}),\rho(T^-_{|_{E^-}})\le 1. \ \ \ \ \ \ \ \ \ \ \  \ \ \ \ \ \ \ \ \ \ \ \ \ \ \ \ \ \ \ \ \ \  \ \ \ \ \  \ $$

\end{defn}

Notice that if $E^- \subset \ho$ is a Banach space such that $S^-(E^-)\subset E$ and $ T^-(E^-)\subset E^-,$ and if the injection $j: E^- \to \ho$ is continuous, then it follows from the closed graph theorem that  the maps $S^-_{|_{E^-}}: E^- \to E^-$ and $ T^-_{|_{E^-}}: E^- \to E^- $ are continuous.

For $M \subset \ho$ set $Z(M):=\{ \zeta \in \C \setminus \overline{\DD} \ | \ f(\zeta)= 0 \ \forall f \in M\},$ and set  ${\mathcal H}_0(\C \setminus \DD):= \cup_{0<r<1}{\mathcal H}_0(\C \setminus r\overline {\DD}),$ so that $\mathcal H(\overline{\DD})\oplus \mathcal H_0 (\C \setminus \DD)={\mathcal O}(\T).$ The argument used in the proof of Proposition 3.2 shows that if $E^- \in \mathcal E^-,$ then $Z(E^-)=\emptyset,$ the injection  $j: E^- \to \ho$ is continuous,  ${\mathcal H}_0(\C \setminus \DD) \subset E^-,$ and $Spec(S_{|_E})=Spec(T_{|_E})=\overline \DD.$

 The following proposition gives an obvious link between the class $\mathcal F$ and the classes $\mathcal E$ and $\mathcal E^-$ of spaces of holomorphic functions introduced in Definition 3.1 and Definition 5.2. 

\begin{prop} (i) Let $F \in \mathcal F.$  The injection $j:F\to \hc$ is continuous, $P^+(F)=F\cap \hd \in \mathcal E$ and $P^-(F)=F\cap \ho \in \mathcal E^-,$ so that ${\mathcal O}(\T) \subset F.$  

Also if $G$ is a closed subspace of $F,$ and if ${\bf S}^{-1}(G) \subset G,$ then $G\cap \hd$ has the division property.

(ii) Conversely if $E \in \mathcal E,$ and if $E^-\in \mathcal E^-,$ then $E\oplus E^- \in \mathcal F.$

\end{prop}

Proof:  

(i) We have, for $f \in \hc,$

$${\bf S}^{-1}P^+f-P^+{\bf S}^{-1}P^+f=\widehat f(0){\bf S}^{-1}1,$$

and so the functional $L_0: f \to \widehat f(0)$ is continuous on $F.$

Hence the functional $L_n: f \to \widehat f(n)$ is continuous on $F$ for $n \in \Z$ since $L_n=L_0\circ S_{|_F}^{-n}.$ This shows that the injection $j: \F \to \hc$ is continuous.

Set $F^+:=F\cap \hd =P^+_{|_F}(F).$ Then $F^+$ is closed.

Let $g \in F^+.$ Then $$Tg={\bf S}^{-1}{\bf S}Tg={\bf S}^{-1}STg={\bf S}^{-1}g -g(0){\bf S}^{-1}1 \in F\cap \hd=F^+.$$

So $T(F^+)\subset F^+,$ and $T_{|_{F^+}}:F^+\to F^+$ is continuous. We have $$\rho(S_{|_{F^+}})\le \rho({\bf S}_{|_F})=1.$$

Let $g \in \hd.$ For $p \ge 0,$ we have $\widehat {T^pg}(n)=0$ for $n<0,$ $\widehat {T^pg}(n)=\widehat g(n+p)=\widehat {{\bf S}^{-p}g}(n)$ for $n\ge 0,$ and so $T^p=P^+\circ ({\bf S}^{-p})_{|_{\hd}}.$
Hence $$\rho(T_{|_{F^+}})\le \rho({\bf S}^{-1}_{|_F})=\rho(({\bf S}_{|_F})^{-1})=1,$$ which shows that $F^+\in \mathcal E.$ So $\mathcal H(\overline {\DD}) \subset F^+.$


Similarly set $F^-:=F\cap \mathcal H_0(\C \setminus \overline{\DD} ) =P^-_{|_F}(F).$ Then $F^-$ is closed, and the same argument as above shows that $F^- \in \mathcal E^-,$ which implies that $F^-$ contains $\mathcal H_0 (\C \setminus \DD).$  Hence $F$ contains $\mathcal H(\overline{\DD})\oplus \mathcal H_0 (\C \setminus \DD)={\mathcal O}(\T).$

Now let $G$ be a closed subspace of $F$ such that ${\bf S}^{-1}(G)\subset G,$ and let $\lambda \in \DD.$ For $g \in G,$ we have

$$({\bf {S}}_{|_F}-\lambda I_{F})^{-1}g=({\bf S}_{|_F})^{-1}(I_F-\lambda ({\bf S}_{|_F})^{-1})^{-1}g=\sum_{n=0}^{+\infty}\lambda^n({\bf S}_{|_F})^{-n-1}g\in G.$$

Now let $g \in G\cap \hd,$ and let $\lambda \in \DD$ such that $g(\lambda)=0.$ We have 

$$g_{\lambda}=({\bf S}_{|_F}-\lambda I_{F})^{-1}({\bf S}_{|_F}-\lambda I_{|_F})g_{\lambda}=({\bf S}_{|_F}-\lambda I_{F})^{-1}g \in G\cap \hd,$$

and so $G\cap \hd$ has the division property. 

\smallskip

(ii) Now let $E\in \mathcal E,$ let $E^- \in \mathcal E^-,$ and set $F:=E\oplus E^-,$ equipped for example with the norm $\Vert f\Vert:=\sqrt{\Vert g \Vert^2+\Vert h\Vert^2}$ for 
$g \in E, h \in E^-.$

Since $\mathcal H(\overline {\DD})\subset E,$ and since $\mathcal H_{0}(\C \setminus \DD) \subset E^-,$ $F$ contains $\mathcal O(\T),$ and in particular ${\bf S}^p.1\in F$ for $p \in \Z.$

 If $g \in E,$ then ${\bf S} g=Sg\in E \subset F,$
and if $h \in E^-,$ then ${\bf S}h= S^-h +\widehat h(-1)1\in F,$ and so ${\bf S}(F) \subset F.$ Similarly if $g \in E,$ then ${\bf S}^{-1}g= Tg +\widehat f(0){\bf S}^{-1}1\in F,$ and if $h\in E^-,$
then ${\bf S}^{-1}h=T^-h\in E^- \subset F,$ and so ${\bf S}^{-1}(F)\subset F.$

Set again $L_n(f)=\widehat f (n)$ for $f \in F, n\in \Z,$ and let $\epsilon >0.$ There exists $C>0$ satisfying $$\max(\Vert (S_{|_E})^n \Vert, \Vert (S^-_{|_{E^-}})^n \Vert, \Vert (T_{|_E})^n\Vert,\Vert (T^-_{|_{E^-}})^n\Vert )\le C(1+\epsilon)^n \ \ (n\ge 0).$$

Let $f \in F,$ and set $g=P^+f \in E,$ $h=P^-f\in E^-.$ We have, for $n\ge 1,$

$$\Vert {\bf S}^ng \Vert =\Vert (S_{|_E})^ng \Vert \le C(1+\epsilon)^n\Vert g\Vert,$$ 
$$ \Vert {\bf S}^nh\Vert \le \Vert (S_{|_{E^-}})^nh\Vert +\sum \limits_{k=1}^n\vert \widehat f (-k)\vert \Vert (S_{|_E})^{n-k}1\Vert$$
$$= \Vert (S_{|_{E^-}})^nh\Vert +\sum \limits_{k=1}^n\vert L_{-1}(({T}_{|_{E^-}}^-)^k) h)\vert \Vert (S_{|_E})^{n-k}1\Vert$$
$$\le C\left (1 +C\Vert {L_{-1}}_{|_{E^-}}\Vert \Vert 1\Vert \right )(1+\epsilon)^n\Vert h\Vert.$$

Hence $\Vert {\bf S}^nf \Vert \le C\left (1 +C\Vert {L_{-1}}_{|_{E^-}}\Vert \Vert 1\Vert \right )(1+\epsilon)^n\Vert f\Vert.$ 

This shows that $\rho \left  ({\bf S}_{|_F}\right )\le 1.$ A similar computation shows that $\rho \left  ({\bf S}^{-1}_{|_F}\right )\le 1.$ So $Spec\left ({\bf S}_F\right )\subset \T,$ and $F\in \mathcal F.$ \ \ $\square$

\smallskip

For $1\le p \le +\infty$ we can identify the space $L^p(\T)$ to a subspace of $\hc$ by identifying $f \in L^p(\T)$ to the hyperfunction $(f^+,f^-)$ where $$f^+(\zeta)=\sum_{n=0}^{+\infty}\widehat f(n)\zeta^n \ \ (\zeta \vert <1), \ \  f^-(\zeta)=\sum_{n=1}^{+\infty}\widehat f(-n)\zeta^{-n}, \ \ (\zeta >1).$$ Since $P^+(L^1(\T))\neq H^1(\DD),$ and since $P^+(L^{\infty}(\T))=BMOA(\DD)\neq H^{\infty}(\DD),$ the spaces $L^1(\T)$ and $L^{\infty}(\T)$ do not belong to the class $\mathcal F.$ 

In the other direction it follows from standard results of harmonic analysis that $L^p(\T)\in \mathcal F$ for $1<p<+\infty.$

Let $\sigma :\Z \to (0,+\infty)$ be a weight satisfying the condition

\begin{equation}0<\inf _{n\in \Z}{\sigma(n+1)\over \sigma(n)}\le \sup_{n\in \Z}{\sigma(n+1)\over \sigma(n)}<+\infty.\end{equation}

Set $\overline \sigma(m)=\sup_{n\in \Z}{\sigma(n)\over \sigma(n+m)}$ and $\tilde \sigma(m)=\sup_{n\in \Z}{\sigma(n+m)\over \sigma(n)}$ for $m\ge 0.$ Then $\overline \sigma$ and $\tilde \sigma$ are submultiplicative, and so the sequences $(\overline \sigma (m)^{1\over m}$ and $(\tilde \sigma (m)^{1\over m}$ have a limit as $m\to +\infty.$ We will denote by 
$\mathcal S$ the class  of weights $\sigma :\Z \to (0,+\infty)$ satifying (6) and the conditions

\begin{equation} \lim_{m\to +\infty}\overline \sigma(m)^{1\over m}=\lim_{m\to +\infty} \tilde \sigma(m)^{1\over m}=1.\end{equation}

For $\sigma \in \mathcal S, 1 \le p < +\infty$ set 
$$HF^p_{\sigma}(\T):= \{ f \in \hc \ | \ \Vert f \Vert_{\sigma, p}:=\left [ \sum \limits_{n\in \Z} \vert \widehat f (n)\vert ^p\sigma(n)^p\right ]^{1\over p}<+\infty \}.$$

A routine well-known verification shows that $HF^p_{\sigma}(\T)\in \mathcal F$ and that we have $$\Vert ({\bf S}_{|_{HF^p_{\sigma}(\T)}})^{-m}\Vert = \overline \sigma(m) \ \mbox{and}
\ \Vert ({\bf S}_{|_{HF^p_{\sigma}(\T)}})^{m}\Vert = \tilde \sigma(m) \ (m\ge 0).$$

Similarly $\mathcal S^+$ denotes the class of weights $\sigma$ on the nonnegative integers satisfying the following conditions

\begin{equation}0<\inf _{n\ge 0}{\sigma(n+1)\over \sigma(n)}\le \sup_{n\ge 0}{\sigma(n+1)\over \sigma(n)}<+\infty.\end{equation}

\begin{equation} \lim_{m\to +\infty}\overline \sigma(m)^{1\over m}=\lim_{m\to +\infty} \tilde \sigma(m)^{1\over m}=1.\end{equation}

where $\overline \sigma(m)=\sup_{n\ge 0}{\sigma(n)\over \sigma(n+m)}$ and $\tilde \sigma(m)=\sup_{n\ge 0}{\sigma(n+m)\over \sigma(n)}$ for $m\ge 0.$

For $1\le p < +\infty,$ we denote by $H^p_{\sigma}(\DD)$ the space of all holomorphic functions $f\to \sum \limits_{n=0}^{+\infty}\widehat f (n)\zeta^n$ holomorphic on the open unit disc
satisfying the condition $\Vert f \Vert_{\sigma,p}:=\left [ \sum \limits_{n=0}^{+\infty}\vert \widehat f (n)\vert^p\sigma(n)^p\right ]^{1\over p}<+\infty.$ Then $H^p_{\sigma}(\DD)\in \mathcal E.$

Also  $\Vert ({T_{|_{H^p_{\sigma}(\DD)}}})^m\Vert = \overline \sigma(m),$ and $\Vert ({S_{|_{H^p_{\sigma}(\DD)}}})^m\Vert =\tilde \sigma(m)$ for $m\ge 0.$

Now denote by $\mathcal S^-$ the class of weights $\sigma$ on the negative integers satisfying the following conditions

\begin{equation}0<\inf _{n \ge 1}{\sigma(-n-1)\over \sigma(-n)}\le \sup_{n\ge 1}{\sigma(-n)\over \sigma(-n-1)}<+\infty.\end{equation}

\begin{equation} \lim_{m\to +\infty}\overline \sigma(m)^{1\over m}=\lim_{m\to +\infty} \tilde \sigma(m)^{1\over m}=1.\end{equation}

where $\overline \sigma(m)=\sup_{n\ge 1}{\sigma(-n)\over \sigma(-n-m)}$ and $\tilde \sigma(m)=\sup_{n\ge 1}{\sigma(-n-m)\over \sigma(-n)}$ for $m\ge 0.$

For $1\le p \le +\infty,$ we denote by $H^p_{0,\sigma}(\C \setminus \overline \DD)$ the space of all holomorphic functions $f\to \sum \limits_{n<0}\widehat f (n)\zeta^n$ holomorphic on $\C\setminus \overline{\DD}$ vanishing at infinity
satisfying the condition $\Vert f \Vert_{\sigma,p}:=\left [ \sum \limits_{n=1}^{+\infty}\vert \widehat f (-n)\vert^p\sigma(-n)^p\right ]^{1\over p}<+\infty.$ 
Then $H^p_{0,\sigma}(\C \setminus \overline \DD)\in \mathcal E^-,$
 $\Vert ({T_{|_{H^p_{0,\sigma}(\C \setminus \overline \DD)}}})^m\Vert = \overline \sigma(m),$ and $\Vert ({S_{|_{H^p_{0,\sigma}(\C \setminus \overline \DD))}}})^m\Vert =\tilde \sigma(m)$ for $m\ge 0.$

The classes $\mathcal S, \mathcal S^+$ and $\mathcal S^-$ are stable under pointwise product, and if $\sigma \in \mathcal S,$ then $(\sigma(n))_{n\ge 0}\in \mathcal S^+$ and 
 $(\sigma(n))_{n<0}\in \mathcal S^-$. Conversely it follows from Proposition 5.3(ii) that if $\sigma$ is a weight on $\Z$ such that $(\sigma(n))_{n\ge 0}\in \mathcal S^+$ and 
 $(\sigma(n))_{n<0}\in \mathcal S^-$, then $\sigma \in \mathcal S,$ which can of course be proved directly by using arguments similar to those used to prove that $E\oplus E^- \in \mathcal F$ if $E \in \mathcal E$ and $E^- \in \mathcal E^-.$ These considerations allow to give to every $E \in \mathcal E$ some Hilbert space "tails", since $E\oplus H^2_{0,\sigma}(\C \setminus \DD)\in \mathcal F$ for every $\sigma \in \mathcal S^-.$
\section{Analytic left-invariant subspaces of Banach spaces of hyperfunctions on the unit circle}

\begin{defn} Let $F\in \mathcal F.$ A linear subspace $N$ of $F$ is said to be left-invariant if ${\bf S}^{-1}(N)\subset N,$ and a left-invariant subspace $N$ of $F$ is said to be analytic if $N\cap\hd\neq 0.$
\end{defn}

Recall that it follows from Proposition 2.4 that if $E \in \mathcal E,$ and if $M$ is a subspace of $E$ having the division property, then there exists for every $\lambda \in \DD$ a linear map $U_{\lambda}: E/M \to E/M$ such that $U_\lambda\circ \pi \circ (S_{|_E}-\lambda I_{E})=\pi,$ where $\pi=E \to E/M$ denotes the canonical surjection. The following result is a general Banach space version of Theorem 3.3 of \cite{ev1}.

\begin{thm} Let $F\in \mathcal F,$ set $L_n(f)=\widehat f(n)$ for $f \in F,$ $n\in \Z$, let $M$ be a closed subspace of $F^+:=F\cap \hd$ having the division property at $0,$ let $\pi:F^+\to F^+/M$ be the canonical surjection, and  let $U_0: F^+/M \to F^+/M$ be the linear map satisfying $U_0\circ \pi \circ S_{|_{F^+}}=\pi.$

If $\sum \limits _{p=1}^{+\infty}\Vert L_{-p}\Vert \Vert U_0^p \pi(1)\Vert<+\infty,$ then $\left [ span \left (\cup_{p\ge 0}{\bf S}_{|_F}^{-p}(M)\right )\right ]^-\cap F^+=M,$ and $\left [span\left ( \cup _{n\ge  0}{\bf S}_{|_F}^{-p}(F^+)\right )\right ]^-=\left [ span \left (\cup_{p\ge0}{\bf S}_{|_F}^{-p}(M)\right )\right ]^-+F^+.$
\end{thm}

Proof: Set $F_1:=\left [ span \left (\cup_{p\ge 0}{\bf S}_{|_F}^{-p}(F^+)\right ) \right ]^-$, and set $N:=\left [  span \left( \cup_{p\ge 0}{\bf S}_{|_F}^{-p}(M)\right ) \right ]^-.$
For $p\ge 1,$ pick $w_p\in \pi^{-1}\left ( U_0^p\pi(1)\right )$ such that $\Vert w_p\Vert \le \Vert U_0^p\pi(1)\Vert +2^{-p}\Vert L_{-p}\Vert^{-1},$
so that $\sum \limits_{p=1}^{+\infty}\Vert L_{-p}\Vert \Vert w_p \Vert<+\infty.$ Define $D :F\to F^+$ by the formula

$$D f=\sum \limits_{p=1}^{+\infty}\widehat f(-p)w_p \ \ (f \in F).$$

Let $u \in M$ such that $u(0)=1,$ and set $v_1=-Tu\in F^+.$ Then $Sv_1=-STu= -u +u(0)1=-u+1,$ and so $1=u+Sv_1.$ 

Set $\alpha_0=1,$ and define by induction $v_p\in M$ for $p\ge 2$ by using the formula $v_p= T(v_{p-1}-v_{p-1}(0)u)$, so that $Sv_p=v_{p-1}-v_{p-1)}(0)u.$ Setting $\alpha_{0}=1, \alpha_p=v_p(0)$ for $p\ge 1, $ we obtain,
\begin{equation} 1= \sum \limits_{k=0}^{p-1}\alpha_k{S}_{|_{F^+}}^ku +S_{|_{F^+}}^pv_p.\end{equation}

An immediate induction shows that $U_0^p\circ \pi \circ S^{p}_{|_{F^+}}=\pi$ for $p\ge 1.$ We obtain

$$U_0^p(\pi(1))= \sum_{k=0}^{p-1}U_0^{p-k}(\pi(u))+\pi(v_p)=\pi(v_p),$$

and so $v_p-w_p\in M$ for $p\ge 1.$ We have

$${\bf S}^{-p}1= \sum_{k=0}^{p-1}\alpha_k{\bf S}^{k-p}u +v_p-w_p+w_p=\sum_{k=0}^{p-1}\alpha_k{\bf S}^{k-p}u +v_p-w_p+D {\bf S}^{-p}1,$$

and so $(P^--D)({\bf S}^{-p}1)\in N$ for $p\ge 1.$

We obtain, for $g\in F^+,$ $p\ge 1,$

$$(P^--D)({\bf S}_{|_F}^{-p}g)=\sum \limits_{k=1}^p\widehat g(k-p)(P^--D)({\bf S}^{-k}1)\in N,$$

and so $(P^--D)(F_1)\subset N,$ and $f =(P^+ +D)f+(P^--D)f\in N+F^+$ for every $f \in F_1.$

To prove the first assertion, consider again $u \in M$ such that $u(0)=1.$ Then $(P^++D)u=u\in M.$ Now assume that $(P^++D){\bf S}_{|_F}^{-k}u\in M$ for $0\le k \le p-1.$ Using (12), we obtain

$$(P^++D){\bf S}_{|_F}^{-p}u=(P^++D)\left ({\bf S}^{-p}1-\sum_{k=1}^{p-1}\alpha_k{\bf S}_{|_F}^{-p+k}u -v_p\right )$$ $$=w_p-v_p-(P^++D)\left (\sum_{k=1}^{p-1}\alpha_k{\bf S}_{|_F}^{-p+k}u \right)
\in w_p-v_p+M\subset M.$$

So $(P^++D){\bf S}_{|_F}^{-p}u\in M$ for $p\ge 1.$ Now let $g \in M.$ Since $M$ has the division property, we see as in the proof of (12) that there exists a sequence $(\beta_p)_{p\ge 0}$ of complex numbers and a sequence $(h_p)_{p\ge 1}$ of elements of $M$ satisfying

$$g=\sum \limits_{k=0}^{p-1}\beta_k{\bf S}_{|_F}^ku+{\bf S}_{|_F}^ph_p \ \ (p\ge 1).$$

Hence $(P^++D){\bf S}_{|_F}^{-p}g \in M$ for $g \in M,$ which implies that $(P^++D)f \in M$ for $f \in \left [ span \left (\cup_{p\ge0}{\bf S}_{|_F}^{-p}(M)\right )\right ]^-.$ Now if $f \in \left [ span \left (\cup_{p\ge0}{\bf S}_{|_F}^{-p}(M)\right )\right ]^-\cap F^+,$ then $f=(P^++D)f\in M,$ and so $M=\left [ span \left (\cup_{p\ge0}{\bf S}_{|_F}^{-p}(M)\right )\right ]^-\cap F^+.$ \ \ $\square$

Notice that if $F \in \mathcal F$ and if $M\subset F^+$ is $z$-invariant and has the division property, then $\left [ span \left (\cup_{p\ge 0}{\bf S}_{|_F}^{-p}(M)\right )\right ]^-$ is a closed
analytic bi-invariant subspace of $F$, which is of course nontrivial if $M $ satisfies the conditions of Theorem 6.2. It follows from a classical result of Wiener that the nontrivial closed biinvariant subspaces of $L^2(\T)$ have the form $N_A:=\{ f \in L^2(\T) \ | \ f (\zeta)=0 \ \mbox{a.e. on A}\}$ for some Lebesgue-measurable subset  $A$ of $\T$ of positive Lebesgue measure, see for example \cite{h}, Chap. 1. Such spaces cannot be analytic since $N_A\cap H^2(\DD)=\emptyset$ if $A$ has positive Lebesgue measure.

In the other direction let $\sigma$ be a nonincreasing unbounded weight on $\Z$ such that $\sigma(n)=1$ for $n\ge 0.$ It follows from the fact that there exist singular inner functions $U$ on the unit disc such that $\vert U^{-1}(z)\vert$ has arbitrarily slow growth as $\vert z \vert \to 1^-$ that the Hilbert space $HF^2_{\sigma}(\T)$ always possessses nontrivial analytic biinvariant subspaces, see \cite{e}.

A natural example of a Banach space $F\in \mathcal F$ satisfying $$F_1:=\left [ span \left (\cup_{p\ge 0}{\bf S}_{|_F}^{-p}(F^+)\right ) \right ]^- \subsetneq F$$ is given by the space $F=H^{\infty}(\DD)\oplus H^{\infty}(\C \setminus \overline{\DD})\subset L^{\infty}(\T).$ In this case $F_1=H^{\infty}(\DD)\oplus \mathcal A_0(\C \setminus \DD),$ where $\mathcal A_0(\C \setminus \DD)$ denotes the algebra of continuous functions on $\C \setminus \DD$ vanishing at infinity which are holomorphic on $\C \setminus \overline{\DD}.$
\section{Consequences of the Matsaiev-Mogulski estimates about the growth of quotient of analytic functions}

In this section, we will use estimates of the growth of analytic functions in the disc which are the quotient of analytic functions satisfying themselves some growth conditions as $\vert \lambda \vert \to 1^-.$
Following the early works of Cartwright \cite{ca} and Linden \cite{l1}, more sophisticated methods were developped in the seventies, see \cite{hk} and \cite{n} for estimates
of inverses of functions analytic in the disc. Concerning analytic functions of the
form $h=f/g,$ where g is allowed to have zeroes, the best results known to the
author are due to Matsaev-Mogulskii \cite{mm}, see also \cite{mo}, who stated their results for functions holomorphic on a half-plane. We state as a theorem the following result, which is the reformulation of Theorem 1 from \cite {mm} given in Corollary 4.2 of \cite{ev1}.
\begin{thm} Let $\Delta: [0,1) \to (0,+\infty)$ be a continuous increasing function and let $f \in \hd.$ Assume that there exists $f_1, f_2\in \hd \setminus \{0\}$ such that $f_2(0)\neq 0, f_2f=f_1$ satisfying for $i=1,2$ the condition

\smallskip

 $$\limsup_{\vert \lambda \vert \to 1^-} \left (log \vert f_i(\lambda)\vert -\Delta(\vert \lambda\vert \right ) <+\infty.$$

\smallskip

(i) If $\int_{0}^1\sqrt{{\Delta(t)\over 1-t}}dt <+\infty,$ then

\begin{equation}log^+\vert f(\lambda)\vert =O\left ({1\over 1-\vert \lambda \vert}\right ) \ \mbox{as} \ \vert \lambda \vert \to 1^-.\end{equation}

\smallskip

(ii) If $\int_{0}^1\sqrt{{\Delta(t)\over 1-t}}dt =+\infty,$ then we have, for every $\epsilon >0,$

\begin{equation} \limsup_{\vert \lambda \vert \to 1^-}\left ( 1-\vert \lambda \vert \right ) log\vert f(\lambda)\vert  \left [ \int_0^{\vert \lambda \vert^{1\over 1+\epsilon}}\sqrt{{\Delta(t)\over 1-t}}dt \right ] ^{-2}\le C(\epsilon),\end{equation}

where \begin{equation}C(\epsilon)={54\over \pi}\epsilon^{-3}(1+\epsilon)\left (1+{2\epsilon\over 3}\right )^2\left (1+{44\over 5}e^{\left (26 \pi +3/2\right )(2+{\epsilon^{-1})}}\right ).\end{equation}

\end{thm}
We now want to give concrete applications of Theorem 6.2 which do not depend on a closed subspace having the division property of a given Banach space $E \in \mathcal E$.

Recall that if $f \in \hd, \lambda \in \DD$ the function $f_\lambda$ is defined for $\zeta \in \DD$ by the formula $f_{\lambda}(\zeta)={f(\zeta)-f(\lambda)\over \zeta -\lambda}$ for $\zeta \neq \lambda$ and $f(\lambda)=f'(\zeta).$ We will use the following easy result.

\begin{lem} Let $E \in \mathcal E,$ let $M$ be a closed subspace of $E$ having the division property, let $\pi: E \to E/M$ be the canonical surjection, let $U:=U_0:E \to E/M$ be the linear map satisfying $U\circ \pi \circ S_{|_E}=\pi,$ and let $f\in M\setminus\{0\}.$ Then we have, for $\lambda \in \DD,$

$$f(\lambda)U(I_{E/M}-\lambda U)^{-1}\pi(1)=-\pi(f_\lambda).$$

\end{lem}

Proof: It follows from corollary 2.6 that for every $\lambda \in \DD$ there exists a linear map $U_{\lambda}:E/M\to E/M$ satisfying $U_{\lambda}\circ \pi\circ (S_{|_E}-\lambda I_E)=\pi,$ that $(I_{E/M}-\lambda U)$ is invertible and that $U_{\lambda}=U(I_{E/M}-\lambda U)^{-1}.$ We obtain

$$f(\lambda)U(I_{E/M}-\lambda U)^{-1}\pi(1)=U_{\lambda}\pi (f(\lambda).1)=U_{\lambda}\pi (f(\lambda).1 -f)$$ $$=-\left (U_{\lambda}\circ \pi \circ (S_{|_E}-\lambda I_E)\right )(f_{\lambda})=-\pi(f_{\lambda}).$$ $\square$

\begin{lem} Let $E \in \mathcal E,$ set $L_0(f) = f(0)$ for $f \in E,$ and set, for $r \in [0, 1),$

$$\Delta_{E}(r):=\sum \limits_{r=0}^{+\infty}r^n\Vert T_{|_E}^n\Vert.$$

Then we have, for $f \in E,$ $\lambda \in \DD,$

$$\vert f(\lambda)\vert \le \Vert L_0\Vert \Vert f \Vert \Delta_{E}(\vert \lambda \vert),  \ \Vert f_{\lambda} \Vert \le \Vert T_{|_E}\Vert \Vert f \Vert \Delta_{E}(\vert \lambda \vert).$$

\end{lem}

Proof: Set $L_n(f)=\widehat f(n)={f^{(n)}(0)\over n!}$ for $f \in E, n\ge 0.$ We have $L_n=L_0\circ T_{|_E}^n,$ so that $\Vert L_n\Vert \le \Vert L_0\Vert \Vert T_{|_E}^n\Vert,$ which gives, for $\lambda \in \DD,$

$$\vert f(\lambda)\vert \le \sum \limits_{n=0}^{+\infty}\vert L_n(f)\vert \lambda \vert ^n= \le \Vert f \Vert \sum \limits_{n=0}^{+\infty}\Vert L_n\Vert \vert \lambda \vert ^n\le \Vert L_0\Vert \Vert f \Vert \Delta_{E}(\vert \lambda \vert).$$

Also it follows from (4) that $f_{\lambda}=T_{|_E}(I_E-\lambda T_{|_E})^{-1},$ which gives, for $\lambda \in \DD,$

$$\Vert f_{\lambda}\Vert = \left \Vert \sum \limits _{n=0}^{+\infty }\lambda ^nT_{|_E}^{n+1}f\right \Vert \le \Vert T_{|_E}\Vert \Vert f \Vert \Delta_E(\vert \lambda \vert).$$ $\square$

Using the Matsaiev-Mogulski estimates, we obtain the following result

\begin{prop}Let $E \in \mathcal E,$ let $M$ be a closed subspace of $E$ having the division property, let $\pi: E \to E/M$ be the canonical surjection, let $U:=U_0:E \to E/M$ be the linear map satisfying $U\circ \pi \circ S_{|_E}=\pi.$ 

(i) If $\int_{0}^1\sqrt {\Delta_E(t)\over 1-t}dt <+\infty,$ then we have, 

\begin{equation}log^{+}\left \Vert \sum \limits _{n=0}^{+\infty}\lambda^nU^{n+1}\pi(1)\right \Vert=O\left ({1\over 1-\vert \lambda \vert }\right ) \ \mbox{as} \ \vert \lambda \vert \to 1^-,\end{equation}

and so

\begin{equation} log^+\Vert U^n\pi(1)\Vert =O\left ({1\over \sqrt n}\right ) \ \mbox{as} \ n \to +\infty.\end{equation}

(ii) If $\int_{0}^1\sqrt {\Delta_E(t)\over 1-t}dt =+\infty,$ then we have, for every $\epsilon >0,$

\begin{equation}\limsup_{\vert \lambda \vert \to 1^-}\left ( 1-\vert \lambda \vert \right ) log \left \Vert \sum \limits _{n=0}^{+\infty}\lambda^nU^{n+1}\pi(1)\right \Vert \left [ \int_0^{\vert \lambda \vert^{1\over 1+\epsilon}}\sqrt{{\Delta_E(t)\over 1-t}}dt \right ] ^{-2}\le C(\epsilon),\end{equation}

where $C(\epsilon)$ is given by (15),  
and so if we set $L_{\epsilon,E}(r) = {C(\epsilon) +1\over 1-r}\left [ \int_0^{r^{1\over 1+\epsilon}}\sqrt{\Delta_E(t)\over 1-t}dt\right ]^2$, we have, when $n$ is sufficiently large,

\begin{equation} \Vert U^n\pi (1)\Vert \le \inf_{0<r<1}r^{-n}e^{L_{\epsilon,E}(r)}.\end{equation}

\end{prop}

Proof: It follows from the Hahn-Banach theorem that the Matsaiev-Mogulski estimates remain valid in the case where the holomorphic functions $f$ and $f_1$ take values in a Banach space $X$, replacing respectively $\vert f (\lambda \vert$ and $\vert f_1 (\lambda \vert)$ by $\Vert f (\lambda \Vert$ and $\Vert f_1 (\lambda \Vert)$,  and (16) and (18)  follow then from Theorem 6.1 and lemmas 6.2 and 6.3.

Now assume that (16) holds. There exists $a>0$ and $b>0$ such that $\left \Vert \sum \limits _{n=0}^{+\infty}\lambda^nU^{n+1}\pi(1)\right \Vert \le ae^{b\over 1-\lambda}$ for $\lambda \in \DD.$ It follows then from the standard vector-valued version of Cauchy's inequalities that we have, for $n\ge 1,$
$$\left \Vert U^n\pi(1)\right \Vert \le \inf_{0< r <1}ar^{-n+1}e^{b\over 1-r}\le a\left (1-{1\over \sqrt n}\right )^{ -n+1}e^{b\sqrt n},$$

and we obtain (17).

Now assume that (18) holds, and let $\epsilon >0.$ It follows from (17) that there exists $r_{\epsilon} \in (0,1)$ such that $log\left \Vert \sum \limits_{n=0}^{+\infty}\lambda^nU^{n+1}\pi(1)\right \Vert
\le L_{\epsilon,E}(\vert \lambda \vert)$ if $ \vert \lambda \vert \in (r_{\epsilon},1).$ There exists $r_n\in (0,1)$ such that $\inf_{0<r<1}r^{-n}e^{L_{\epsilon}(r)}=r_n^{-n}e^{L_{\epsilon,E}(r_n)},$
and $\lim_{n\to +\infty}r_n=1$ since $L(r) \ge {r^n\over r_n^n}L(r_n)$ for $r_n < r < 1.$ Let $p\ge 1$ such that $r_n >r_{\epsilon}$ for $n \ge p.$ It follows from the vector valued version of Cauchy's inequalities that we have

$$\Vert U^n\pi(1)\Vert \le r_n^{-n+1}\sup_{\vert \lambda \vert =r_n}\left \Vert \sum \limits_{m=0}^{+\infty}\lambda^mU^{m+1}\pi(1)\right \Vert\le r_n^{-n}e^{L_{\epsilon,E}(r_n)} \ \ (n\ge p).$$

  $\square$

Recall that a sequence $(u_n)_{n\ge p}$ of positive real numbers is said to be {\it log-concave} if the sequence $({u_{n+1}\over u_n})_{n\ge p}$ is nonincreasing, and that a sequence 
$(u_n)_{n\ge 1}$ is said to be {\it eventually log-concave}  is the sequence $(u_n)_{n\ge p}$  is log-concave for some $p\ge 1.$ 

Let $(u_n)_{n\ge 0}$ be an eventually log-concave sequence $(u_n)_{n\ge 0}$ of real numbers, satisfying $u_n\ge 1$ for $n\ge 0,$ satisfying the one-sided "nonquasianalytic condition"

\begin{equation} \sum \limits_{n=1}^{+\infty} {log (u_n)\over n^{3\over 2}}<+\infty,\end{equation}

 set  $v_n= (n+1)^2u_n$ for $n\ge 0,$ and set $\Lambda(r)=\sup_{n\ge 0}r^nv_n$ for $r \in [0,1[$. Since the sequence $(v_n)_{n\ge 0}$ is also eventually log-concave  and satisfies (20), it follows from \cite{n}, section 2.6, lemma 2 that $\int_0^1\sqrt {log \Lambda(t)\over 1-t}dt<+\infty.$ 
 
 So if $E \in \mathcal E,$ and if $\Vert T_{|_E}^n\Vert =O(u_n)$ as $n\to +\infty$, there exists $k>0$ such that we have, for $r \in [0,1),$
 
 $$\Delta_E(r)=\sum \limits_{n=0}^{+\infty}r^n\Vert T_E^n\Vert \le \left ( \sum_{n=0}^{+\infty}{1\over (n+1)^2}\right )\sup_{n\ge 0}r^n(n+1)^2\Vert T_{|_E}^n\Vert\le k{\pi^2\over 6}\Lambda(r),$$

and so $\int_0^1\sqrt{\Delta_E(t)\over 1-t}dt<+\infty,$ and $log^+\Vert U^n\pi(1)\Vert=O\left ({1\over \sqrt n}\right )$ for every nontrivial closed subspace $M$ of $E$ having the division property.

If $E \in \mathcal E,$ and if $\int_0^1\sqrt{\Delta_E(t)\over 1-t}dt=+\infty,$ set $\sigma(n)=\sigma_{\epsilon, E}(n):=\inf_{0<r<1}r^{-n}e^{L_{\epsilon, E}(r)}$ for $n\ge 0.$ Then the sequence $\sigma$ is increasing and log-concave, $\sigma(0)=1,$ and ${\sigma(n+p)\over \sigma(p)}={\sigma(n+p)\over \sigma(n+p-1)}\dots {\sigma(n+1)\over \sigma(n)}\le {\sigma(p)\over \sigma(p-1)}\dots{\sigma(1)\over \sigma(0)}=\sigma(p),$ so that, with the notations of (8) and (9), we have $\overline \sigma (n)=1$, $\tilde \sigma (n)=\sigma(n)$ for $n\ge 0,$ and, clearly, $\lim_{n\to +\infty}\sigma(n)^{1\over n}=1.$ So, with the notations of section 5, we have $\sigma \in \mathcal S^+.$ 

Now set $\check{\sigma}_{\epsilon, E}(-n)=(n+1)^2\sigma (n)$ for $n\ge 1,$ and for $1\le p <+\infty$ consider the weighted Hardy space $H^p_{0,\check {\sigma}_{\epsilon, E}}(\C \setminus \overline{\DD})$ introduced in section 5. It follows from Proposition 7.4 that if $M$ is a closed nontrivial subspace of $E$ having the division property, then 

$$\sum \limits_{n=1}^{+\infty}\Vert L_{-n}\Vert\Vert U^n\pi(1)\Vert=\sum \limits_{n=1}^{+\infty}\check{\sigma}(-n)^{-1}\Vert U^n\pi(1)\Vert<+\infty,$$

where $L_n(f)=\widehat f(n)$ for $f \in H^p_{0,\check {\sigma}_{\epsilon, E}}(\C \setminus \overline{\DD}),$ $n\le -1$,
and it follows from Theorem 6.2 that if we set $F:=E\oplus H^p_{0,\check {\sigma}_{\epsilon, E}}(\C \setminus \overline{\DD}),$ then  we have 

\begin{equation}F=\left [ span \left (\cup_{p\ge0}{\bf S}_{|_F}^{-p}(M)\right )\right ]^-+F^+, \left [ span \left (\cup_{p\ge 0}{\bf S}_{|_F}^{-p}(M)\right )\right ]^-\cap F^+=M.\end{equation}
 
 A tedious verification that we omit shows that this result remains true if we just set $\check{\sigma}_{\epsilon, E}(-n)=\sigma_{\epsilon, E} (n)$ for $n\ge 1.$

For example assume that $E \in \mathcal E,$ and that $log\Vert T_{|_E}^n\Vert =O(n^\alpha)$ as $n\to +\infty,$ where $\alpha \in [0,1).$ Let $\sigma \in \mathcal S^-.$ It follows from the above discussion and from computations given in \cite{ev1}, section 4, that every closed subspace of $E$ having the division property satisfies (21) with respect to $F=E\oplus H^p_{\sigma}(\C \setminus \overline{\DD}), 1 \le p <+\infty$ if $\sigma$ satisfies the following conditions

\smallskip

- $\liminf_{n\to +\infty}{log (\sigma(-n))\over \sqrt n}= +\infty$ when $0 \le \alpha <{1\over 2}$

- $\liminf_{n\to +\infty}{log (\sigma(-n))\over \sqrt n log(n+1)}= +\infty$ when $\alpha ={1\over 2}$

- $\liminf_{n\to +\infty}{log (\sigma(-n))\over n^{\alpha}}= +\infty$ when ${1\over 2}< \alpha <1.$

Now assume that $log\Vert T_{|_E}^n\Vert= O\left (n\over (log(n+1))^c\right )$ as $n\to +\infty$ for some $c>0.$ 

It follows also from the above discussion and from computations given in \cite{ev1}, section 4 that every closed subspace of $E$ having the division property satisfies (21) with respect to $F=E\oplus H^p_{\sigma}(\C \setminus \overline{\DD}), 1 \le p <+\infty$ if  $\sigma \in \mathcal S^-$ satisfies the condition

$$\liminf_{n \to +\infty}{log (\sigma(-n)(log(n+1))^c\over n}>\limsup_{n\to +\infty}{log \Vert T_{|_E}^n\Vert(log(n+1))^c\over n}.$$

We conclude this section with the following other consequence of the Matsaiev-Mogulski estimates, which will be used in the next section.

\begin{prop} Let $E \in \mathcal E,$ let $\sigma \in \mathcal S^-,$ and set $F=E\oplus H^p_{0,\sigma}(\C \setminus \overline{\DD}),$ where $p\ge 1.$  Assume that $(\sigma(-n))_{n\ge 1}$ is eventually log-concave, and that the sequence $\left ( {log(\sigma(-n))\over n^{\alpha} }\right )_{n\ge 1}$ is eventually nondecreasing for every $\alpha > 0.$

If $\limsup_{n\to +\infty} {log\Vert T_E^n\Vert \over log(\sigma(-n))}<1,$ then $N=\left [ \cup_{n\ge 0} {\bf S}^{-n}(N\cap E)\right ]^-$ and $F=N+E$ for evry analytic closed left-invariant subspace $N$ of $F.$

\end{prop}

Proof: There exists $c\in (0,1)$ and $m>0$ such that $\Delta_E(r)\le m\sum \limits_{n=0}\sigma^c(-n-1)$ for $r\in (0,1).$ 

Using (18) and (19), Cauchy's inequalities and estimates given in the proof of Proposition 4.4 of \cite{ev1} concerning the growth of ${1\over 1-r}\int_{0}^{r^{1\over {1+\epsilon}}}\sqrt{ {\sum \limits_{n=0}^{+\infty}t^n\sigma^c(-n-1)}\over 1-t}dt,$ we see that $\sum \limits_{n=1}^{+\infty}\Vert U_M^n\Vert \sigma(-n)<+\infty,$ for every nontrivial closed subspace of $E$ having the division property, and the result follows from Theorem  6.2. $\square$

\section{Closed $z$-invariant subspaces having the division property and closed biinvariant subspaces of Banach spaces of hyperfunctions on the unit circle}

In this final section we indicate how to associate to every Banach space of holomorphic functions $E \in \mathcal E$ a weight $\sigma \in \mathcal S^-$ such that for every nontrivial closed  subspace $N$ of $E\oplus H^2_{0,\sigma}$ which is invariant for $\bf S^{-1}$ there exists $k \ge 0$ such that ${\bf S}^k(N)$ is analytic. This means that  ${\bf S}^k(N)\cap E\neq \{0\},$ for some $k \ge 0,$ and we have more precisely in this case

$$N= \left [span \left (\cup_{p\ge 0}{\bf S}_{|_F}^{-p}(({\bf S}^k_{|_F}(N)\cap E)\right )\right ]^-, $$ $$E \oplus H^2_{0,\sigma}=\left [ span \left (\cup_{p\ge 0}{\bf S}^{-p}_{|_F}({\bf S}^k_{|_F}(N\cap E)\right )\right ]^-+E.$$

In particular if $N$ is a closed nontrivial invariant subspace of $E\oplus H^2_{0,\sigma}$ which is invariant for $\bf S$ and $\bf S^{-1},$ then $N= \left [span \left (\cup_{p\ge 0}{\bf S}_{|_F}^{-p}(N\cap E)\right )\right ]^-.$

This result was proved in \cite{ev1} when $E=H^2_{\tau}(\DD),$ where $\tau \in {\mathcal S}^+$, and the proof relies heavily on the theory of asymptotically holomorphic functions in the disc. We will give here the modifications needed to extend the results of \cite{ev1} to all Banach spaces $E\in \mathcal E.$

If $F \in \mathcal F,$ and if $w \in \hc$ satisfies the condition $\sum \limits_{n\in \Z}\vert \widehat w(n)\vert \Vert {\bf S}^n_{|_F}\Vert <+\infty,$ set $w({\bf S}_{|_F}):=\sum \limits _{n\in \Z}\widehat w(n){\bf S}^n_{|_F}.$ We obtain, for $f \in F,$

\begin{equation} \widehat {w({\bf S}_{|_F})f}=\widehat w *\widehat f \ \ (f\in F),\end{equation}
where the coefficients of the convolution product are given by absolutely convergent series. Notice that (22)  implies that  $w=w({\bf S}_{|_F})1\in F.$ This allows in particular to define $w({\bf S}_{|_F})$ for every $w \in {\mathcal O}(\T)$ and every $F\in \mathcal F.$

 In order to present the strategy of the proof, we first give a consequence of a standard factorization result for functions holomorphic in an annulus, which was given in \cite{ev1}, Prop. B1 in the case where $F=HF^2_{\sigma}(\T)$ with $\sigma \in \mathcal S.$

\begin{prop} Let $F \in \mathcal F,$ let $f =(f^+,f^-)\in F,$ and assume that $f^-\in \mathcal H(\C\setminus {{\DD}}).$  Then there exists a function $g \in F^+,$ a nonnegative integer $k$  and a function $h \in \mathcal H_0(\C\setminus {{\DD}})$ such that $f=e^{h({\bf S_{|_F}^{-1}})}{{\bf S}_{|_F}^{-k}}g,$ and we have

\begin{equation}\left [ span\left \{ {\bf S}_{|_F}^{-n-k}g\}_{n\ge 0}\right \} \right ]^-=\left [ span\left \{ {\bf S}_{|_F}^{-n}f\}_{n\ge 0}\right \} \right ]^-.\end{equation}

\end{prop}

Proof: There exist $r_0\in (0,1)$ such that $f^-$ has a holomorphic extension to $\C \setminus r_0\overline {\DD},$ and so we can set $f(\zeta)=f^+(\zeta)+f^-(\zeta)$ for $z\in \DD \setminus r_0\overline{\DD}.$ Pick $r \in (r_0,1).$ It follows from standard complex analysis results that there exists $g \in \hd, h \in {\mathcal H}(\C \setminus r{\overline \DD})$ and a nonnegative integer $k$ such that $f(\zeta)=\zeta^{-k}e^{h(\zeta)}g(\zeta)$ for $\zeta \in \DD \setminus r{\overline \DD},$ so that $g(\zeta)= e^{-h(\zeta)}\zeta^kf(\zeta).$ 

Set $u(\zeta)=\zeta^kf(\zeta), v(\zeta)=e^{-h(\zeta)} $ for $\zeta \in \DD \setminus r\overline {\DD}.$ The Fourier coefficients of elements of $\mathcal H(\DD \setminus r\overline{\DD}),$ considered as elements of $\hc$ coincide with the coefficients of their Laurent series expansion, and so we have $\widehat g=\widehat v *\widehat u.$ It follows then from (22) that $g=( e^{-h})({\bf S}_{|_F}){\bf S}_{|_F}^pf\in F^+.$ 

Property (23) follows then from the fact that  $g=(e^{-h})({\bf S}_{|_F)}){\bf S}_{|_F}^kf$ and $f =(e^{h})({\bf S}_{|_F)}){\bf S}_{|_F}^{-k}g.$ \ \ $\square$

\smallskip

An obvious idea to construct weights $\sigma\in \mathcal S^-$ such that   there exists $k\ge 0$ and $g \in E$ satisfying (23) for $f \in F:=E\oplus H^2_{0,\sigma}$ would be to obtain for every $f \in F$ factorizations of the form $f=e^{h({\bf S}_{|_F})}{\bf S}^{-k}g,$ where $g \in E,$ where $\widehat h(n)=0$ for $n\ge 0,$ and where $h({\bf S}_{|_F})$ is a bounded operator on $F.$ Unfortunately, the author showed in \cite{e1}, using basic facts about Banach algebras, that such nice factorizations cannot hold for every $f\in F.$ So as in \cite{ev1} we will indeed obtain when $\sigma$ grows "sufficiently rapidly and regularly" factorizations of the same type, but where $h({\bf S}_{|_F})$ maps $F$ into a larger Banach space $\tilde F=E\oplus H^2_{0,\tilde \sigma}$, where $\lim _{n\to -\infty}{\tilde \sigma(n)\over \sigma(n)} = 0.$ 

The following lemma is an extension to all Banach spaces $E \in \mathcal E$ of a result given in \cite{ev1}, Prop. 5.1 for weighted Hardy spaces $H^2_{\tau}(\DD).$ where $\tau \in \mathcal S^+.$

\begin{lem} Let $E \in \mathcal E,$ let $\sigma\in \mathcal E^-,$ and set $F:=E\oplus H^2_{0,\sigma}$. Assume that $\left (\sigma(-n)\right )_{n\ge 1}$ is eventually log-concave and satisfies the following conditions

(i) the sequence $\left ({log(\sigma(-n)\over n^{a}}\right )_{n \ge 1}$ is eventually nondecreasing for every $a \in (0,1),$

(ii) $ \limsup_{n\to +\infty}{log \left \Vert T_{|_F}^n\right \Vert\over log \left(\sigma(-n)\right )}\le 1.$

Then we have

$$\limsup_{n\to +\infty}{log \left \Vert {\bf S}^{-n}_{|_F}\right \Vert \over log \left(\sigma(-n)\right )} = 1.$$

\end{lem}

Proof: Equip for example $F$ with the norm $\Vert(f^+, f^-)\Vert =\sqrt{\Vert f^+\Vert^2+\Vert f^-\Vert^2}.$ We have, for $f=(f^+,f^-) \in F,$ $n \ge 1,$
 
$${\bf S}^{-n}f=(Tf^+,\sum \limits_{m=0}^{n-1}\widehat f(m){\bf S}^{m-n}1 +{\bf S}^{-1}f^-),$$

$$\Vert {\bf S}^{-n}f\Vert^2=\Vert Tf^+\Vert ^2 +\sum \limits_{m=0}^{n-1}\vert \widehat f(m)\vert^2\sigma^2(m-n) +\Vert {\bf S}^{-1}f^-\Vert^2.$$
Set, for $g\in E,$

 $$R_n(g):=\sum \limits_{m=0}^{n-1}\widehat g(m){\bf S}^{m-n}1\in H^2_{0,\sigma}.$$
 
 We obtain
 
 $$\Vert {\bf S}^{-n}_{|_F}\Vert =\max \left (\left \Vert T^n_{|_E}\right \Vert, \left \Vert R_n\right \Vert, \left \Vert {\bf S}^{-n}_{|_{H^2_{0,\sigma}}}\right \Vert\right ).$$
 
 We have $\left \Vert {\bf S}^{-n}_{|_{H^2_{0,\sigma}}}\right \Vert=\sup_{p\ge 0}{\sigma(-n-p)\over \sigma(-p)}.$ Since the sequence $((\sigma(-n))_{n\ge 1}$ is eventually log-concave, an elementary well-known computation shows that there exists two positive reals $\alpha, \beta$ such that $\alpha \sigma(-n)\le \sup_{p\ge 0}{\sigma(-n-p)\over \sigma(-p)}\le \beta \sigma(-n)$ for $n\ge 1,$ and so $\lim_{n\to +\infty} {log \left \Vert {\bf S}^{-n}_{|_{H^2_{0,\sigma}}}\right \Vert\over log (\sigma(-n))}=1.$
 
 It remains to show that $\limsup_{n\to +\infty}{log \Vert R_n\Vert\over log(\sigma(-n))}\le 1.$ Let $\epsilon >0,$ and let $a \in (0,1)$ such that $\left (1+{\epsilon\over 2}\right )2^{1-a} <1+ \epsilon.$ Changing a finite numbers of terms of the sequence $(\sigma(-n))_{n\ge 1},$ which does not affect $\limsup_{n\to +\infty}{log \left \Vert {\bf S}^{-n}_{|_F}\right \Vert \over log \left(\sigma(-n)\right )},$ we may assume that the sequence $\left ({log (\sigma(-n))\over n^a}\right)_{n\ge 1}$ is nondecreasing, and there exists 
 $\mu >0$ such that $\Vert T_{|_E}^n\Vert \le \mu \sigma(-n)^{1+{\epsilon\over 2}}$ for $n\ge 1.$ Since $x^a +(1-x)^a\le 2 ^{1-a}$ for $0\le x \le 1,$ we obtain, for $0 \le m \le n-1,$ $$log(\sigma(-m)) +log (\sigma(n-m))\le \left ({m^{a} +(n-m)^a\over n^a}\right ) log(\sigma(-n))\le 2^{1-a}log(\sigma(-n)).$$
 
 Set $L_m(g)=\widehat g(m)$ for $f \in E, m\ge 0.$ Since $L_m=L_0\circ T_{|_E}^m,$ we have
 
 $$\Vert R_n \Vert \le \sum _{m=0}^{n-1}\Vert L_0\Vert \Vert T^m_{|_E}\Vert \sigma(n-m)\le n \mu \Vert L_0\Vert \sup_{0\le m \le n-1}\sigma(-m)^{1+{\epsilon \over 2}} \sigma(n-m)^{1+{\epsilon\over 2}},$$
 
 and so $\limsup_{n\to +\infty}{log \Vert R_n\Vert\over log(\sigma(-n))}< 1+\epsilon.$ $\square$

We now give an extension to Banach spaces $E \in \mathcal E$ of  a result obtained in \cite{ev1}, corollary 5.3 for weighted $H^2$-spaces of holomorphic functions in the open unit disc. A seminal result in this direction was given in \cite{bv}, Theorem 6.1.
 
 \begin{thm} Let $E\in \mathcal E,$ and assume that a weight $\sigma \in \mathcal S^-$ satisfies the following conditions
 
 (i) the sequence $(\sigma(-n))_{n\ge 1}$ is eventually log-concave,
 
 (ii) the sequence $\left ({\log (\sigma(-n))\over n}log(n)^A\right )_{n\ge 1}$ is eventually nondecreasing for some $A>0,$
 
 (iii) $\sum \limits_{n<0}{log(\sigma(-n))\over n^2}=+\infty,$
 
 (iv) $\limsup_{n\to +\infty}{log \Vert T_{|_E}^n\Vert \over log(\sigma(-n))}<{1\over 64}.$
 
 Set $F_s:=E\oplus H^2_{0, \sigma^s}$ for $s>0,$ $F:=F_1.$ Then for every $f \in F$ and  for every $s <{1\over 4}$ there exist $k\ge 0$, $g \in E$ and a function $w \in \ho$ satisfying the following conditions
 
 (1) $\sum \limits _{n<0}\vert \widehat w (n)\vert \Vert {\bf S}^n_{|_{F_s}}\Vert <+\infty,$

 (2) $f=e^{w \left ({\bf S}_{|_{F_s}}\right )}{\bf S}^{-k}_{|_{F_s}}g.$

 \end{thm}
 Proof: Set $\sigma(n)={1\over (n+1)\Vert T_{|_E}^n\Vert}$ for $n \ge 0.$ Then $$\overline {\sigma}^+(n):=\sup_{p\ge 0}{\sigma(p)\over \sigma(n+p)}=\sup_{p\ge 0}{(n+p+1)\Vert T_{|_E}^{n+p}\Vert \over (p+1)\Vert T_{|_E}^p\Vert}\le (n+1)\Vert T_{|_E}^n\Vert \ \ (n\ge 0),$$
 and so $\limsup_{n\to +\infty}{log \overline{\sigma}^+(n)\Vert \over log(\sigma(-n))}<{1\over 64}.$ The weight $\sigma$ satisfies the conditions of Corollary 5.3 of \cite{ev1}, and so for every $s<{1\over 4}$ and every $f \in HF^2_{\sigma}$ there exists $w\in \ho$, $k \ge 0$ and $g\in H^2_{\sigma^+}(\DD)$ such that $\sum_{n<0}\vert \widehat w(n)\vert^2 \sigma^{2s}(n)<+\infty$ satisfying the condition $\widehat f =\widehat {e^w}*\widehat {{\bf S}^{-k}g}.$
 
 By construction, we have $E\subset H_{\sigma}^2(\DD).$  Let $s \in (0, {1\over 4}),$ choose $s_0 \in (s, {1\over 4}),$ and let $f \in HF^2_{\sigma}(\T).$ Apply the factorization result at $f$ and $s_0. $ and set $F_s:=E\oplus H^2_{0, \sigma^s}.$ It  follows from the lemma that we have
 
 $$\limsup_{n\to +\infty}{log\Vert {\bf S}^{-n}_{|_{F_s}}\Vert\over log(\sigma^{s_0}(-n))}={s_0\over s}\limsup_{n\to +\infty}{log\Vert {\bf S}^{-n}_{|_{F_s}}\Vert\over log(\sigma^{s}(-n))}\le {s_0\over s}<1.$$
 
 So $\sum_{n<0}\widehat w(n) \Vert {\bf S}^n_{|_{F_s}}\Vert<+\infty,$  and it follows from (22) that $f =e^{w({\bf S}_{|_{F_s}})}{\bf S}^k_{|_{F_s}}g,$ and
  $g=e^{-w({\bf S}_{|_{F_s}})}{\bf S}^k_{|_{F_s}}f\in E.$ $\square$

  In the following if $F=E\oplus H^2_{0,\sigma}(\C \setminus \overline {\DD}),$ where $E\in \mathcal E$ is given, and where $\sigma \in \mathcal S^-$ is eventually log-concave, the dual space $F^*$ is identified to the space $E^*\oplus H^2_{\sigma^*}(\DD),$ where $E^*$ denotes the dual space of $E,$ where $\sigma^*(n)={1\over \sigma(-n-1)}$ for $n\ge 0,$ and where $H^2_{\sigma^*}(\DD)$ is identified with the dual of $H^2_{0,\sigma}(\C \setminus \overline {\DD}),$ the duality being implemented by the formula
  
  $$<u, v>:=\sum \limits_{n=0}^{+\infty}\widehat u(-n-1)\widehat v (n) \ \ (u \in H^2_{0,\sigma}(\C \setminus \overline{\DD}), v \in H^2_{\sigma^*}(\DD)).$$
  
  We now present the so-called Dynkin transform, that we will use as in \cite{ev1}. A weight $\tau$ on $\Z^+$ is said to be  log-convex when $\tau^{-1}$ is log-concave, and two weights $\tau_1$ and $\tau_2$ on $\Z^+$ are said to be equivalent when $$0 < \inf_{n \in \Z^+}{\tau_1(n)\over \tau_2(n)}\le \sup_{n \in \Z^+}{\tau_1(n)\over \tau_2(n)}<+\infty.$$

  Set $L^2_+([0,1]):=\{\phi \in L^2([0,1]) \ | \ \phi(t)>0 \ a.e.\},$ and for $\phi \in L^2_+([0,1]),$ set
 
  \begin{equation} \tau_{\phi}(n)=\left [ 2\int_0^1\phi^2(t)t^{2n+1}dt\right ]^{1\over 2} \ \ (n\ge 0).\end{equation}

  Then $\tau_{\phi}$ is log-convex, and conversely, if a weight $\tau$ on $\Z^+$ is eventually log-convex,  it follows from Appendix A of \cite{b} and from  lemma 5.2 of \cite{d} that there exists $\phi \in L^2_+([0,1])$
  such that $\tau_{\phi}$ is equivalent to $\tau.$ So if $(\sigma(-n))_{n\ge 1}$ is eventually log-concave, the set $W({\sigma})$ of all functions $\phi \in L^2_+([0,1])$ such that $\tau_{\phi}$ is equivalent to $\sigma^*$ is nonempty.
  
The usual Cauchy transform of $F \in L^1(\DD)$ is defined by the formula

$$\mathcal C(F)(\lambda):={1\over \pi}\int \int _{\DD}{F(\zeta) \over \lambda-\zeta}dm(\zeta),$$

where $m$ denotes the $2$-dimensional Lebesgue measure. 

We have $\mathcal C(F) \in L^p_{loc}(\C)$ for $1\le p \le 2,$ and $\mathcal C (F)$ is bounded and continuous on $\C$ if, further, $F\in \cup_{q>2}L^q(\C).$ Also if  $\mathcal C^+(F)$ and  $\mathcal C^-(F)$ denote the restrictions of $\mathcal C(F)$ to $\DD$ and $\C \setminus \overline {\DD},$ then
 $\overline \partial \mathcal C^+(F)=F$ in the sense of distribution theory, and $\mathcal C^-(F)\in \ho.$
 
 Let $\phi \in W(\sigma),$ set $\mathcal B_{\phi}:=\{ F \in \hd \ | \ \int \int_{\DD} \vert F(\zeta)\vert^2\phi^2(\vert \zeta \vert)dm(\zeta)<+\infty \},$ set $\phi(\zeta) =\phi(\vert \zeta \vert)$ for $\zeta \in \DD,$ and set
 
 $$ [F,G]:={1\over \pi}\int \int_{\DD} F(\zeta)\overline G(\zeta)dm(\zeta) \ \  (F\in \mathcal B_{\phi}, G \in \phi^2\overline{\mathcal B_{\phi}}).$$
 
 This allows to identify $\phi^2\overline {\mathcal B_{\phi}}$ to the dual space of $\mathcal B_{\phi},$ and $\phi^2\overline {{\mathcal B}_{\phi}}\subset L^1(\DD)$ since $\phi \in L^2(\DD)$ and $\phi \overline F\in L^2(\DD)$ for $F \in \mathcal B_{\phi}.$
 
Now let $\phi \in W(\sigma).$ Then $\mathcal B_{\phi}=H^2_{\sigma*}(\DD),$ and the $\phi$-Dynkin transdorm $\mathcal D_{\phi}(h)$ is defined for $h \in H^2_{0,\sigma}(\overline \DD)$ by the formula

\begin{equation} \mathcal D_{\phi}(h)=\mathcal C^+(L_{\phi}(h)),\end{equation},

where $L_{\phi}(h)$ is the unique element of $\phi^2\overline{\mathcal B_{\phi}}$ satisfying

\begin{equation}  [F,L_{\phi}(h)]=\langle h, F\rangle \ \ \forall F \in \mathcal B_{\phi}=H^2_{\sigma^*}(\DD).\end{equation}

Then $h=\mathcal C^-(L_{\phi}(h))$, and this  is why $\mathcal D_{\phi}(h)$ is called an extension of $h$ to $\DD.$

The following result is an extension of Theorem 3.1 of \cite{ev1} to the case where $E\in \mathcal E.$ 

\begin{prop} Let $E \in \mathcal E,$ let $\sigma \in \mathcal S^-$ be eventually log-concave, let $f =(g,h)\in E\oplus H^2_{0,\sigma}(\C \setminus \overline{\DD}),$ let $l=(v, \theta)\in E^*\oplus H^2_{\sigma*}(\DD),$ and let $\phi \in W(\sigma).$

Then the following conditions imply each other

\smallskip
(i) $\langle{\bf S}^{-n}f,l\rangle=0 \ \ \mbox{for} \  n\ge 1,$

\smallskip
(ii) $\theta(g+\mathcal D_{\phi}(h))=\mathcal C^+(\theta\overline \partial \mathcal D_{\phi}(h))-\langle g_.,v\rangle ,$

where $\langle g_.,v \rangle(\lambda)=\langle g_{\lambda},v \rangle $ for $\lambda \in \DD.$

\end{prop}

Proof: For the convenience of the reader, we give a self-contained proof, which seems somewhat simpler than the proof given in \cite{ev1} when $E$ is a weighted $H^2$-space of holomorphic functions in the unit disc. Since $f=(g,h),$ we have, for $n\ge 1,$

$${\bf S}^{-n}f=\left (T^ng, R_ng +{\bf S}^{-n}h\right),$$

where $R_n:E \to H^2_{0,\sigma}$  is defined for $u \in E$ by the formula

$$R_nu =\sum_{m=0}^{n-1}\widehat u(m){\bf S}^{m-n}.1,$$

so that $\langle \theta, R_nu \rangle=\sum \limits_{m=0}^{n-1}\widehat u(m)\widehat \theta(n-m-1)=\widehat{u\theta}(n-1).$

So the series $\sum \limits _{n=0}^{+\infty}\lambda^{n}\langle \theta, R_{n+1}g \rangle$ converges for $\lambda \in \DD,$ and we have

$$\sum \limits _{n=0}^{+\infty}\lambda^{n}\langle \theta, R_{n+1}g \rangle=g(\lambda)\theta(\lambda) \ \ (\lambda \in \DD).$$

The series $\sum \limits _{n=0}^{+\infty}\lambda^{n}\langle T_{|_E}^{n+1}g,v \rangle =\left \langle \sum \limits_{n=0}^{+\infty}\lambda^{n}T_{|_E}^{n+1}g, v \right \rangle$ converges for every $\lambda \in \DD,$ and it follows from (4) that we have

$$\sum \limits _{n=0}^{+\infty}\lambda^{n}\langle T_{|_E}^{n+1}g,v \rangle=\left \langle \sum \limits _{n=0}^{+\infty}\lambda^{n}T_{|_E}^{n+1}g,v \right \rangle=\langle g_{\lambda}, v \rangle \ \ (\lambda \in \DD).$$

Also, with the convention $\widehat h (p)=0$ for $p\ge 0,$ we have, for $n\ge 1,$

$$\langle  \theta, {\bf S}^{-n}h \rangle= \sum \limits  _{m=n}^{+\infty}\widehat \theta (m)\widehat h(n-m-1)
=\sum \limits_{p=0}^{+\infty}\widehat \theta (p+n)\widehat h(-p-1)=\langle T^n\theta,h \rangle.$$

We obtain, for $\lambda \in \DD,$

$$\sum \limits _{n=0}^{+\infty}\lambda^n\left \langle  \theta, S^{-n-1}u\right \rangle=\left \langle \sum \limits_{n=0}^{+\infty} \lambda^n T^{n+1},h\right \rangle=<\theta_{\lambda},h> $$ $$={1\over \pi}\int \int_{\DD} {\theta (\zeta)-\theta(\lambda)\over \zeta -\lambda}\overline \partial\mathcal D_{\phi}(h)dm(\zeta)=\theta(\lambda)\mathcal D_{\phi}(h)(\lambda)-\mathcal C^+(\theta \overline \partial \mathcal D_{\phi}(h))(\lambda),$$

Since $<{\bf S}^{-n}f,l>= <T^ng,v>+ <\theta, R_ng>+ <\theta, {\bf S}^{-n}h>$ for $n \ge 1,$ the result follows.\ \
 $\square$
 
 Recall that a linear space $N \subset \hc$ is said to be left-invariant if ${\bf S}^{-1}(N)\subset N,$ and $N$ is said to be bi-invariant if ${\bf S}(N)\cup {\bf S}^{-1}(N)\subset N.$ Also a left-invariant subpace $N$ of a Banach space $F\in \mathcal F$ is said to be 
 analytic if $N\cap \hd=N\cap F^+\neq \{0\}.$
 We are now ready to associate to every Banach space $E\in \mathcal E$ a Hilbert space "tail" $H^2_{\sigma}(\C \setminus \overline \DD)$ for which every nontrivial bi-invariant closed subspace of $E\oplus H^2_{\sigma}(\C \setminus \overline \DD)$ is generated by a nontrivial closed $z$-invariant subspace of $E$ having the division property.

 \begin{thm} Let $E\in \mathcal E,$ and assume that a weight $\sigma \in \mathcal S^-$ satisfies the following conditions
 
 (i) the sequence $({\sigma(-n)\over n^{\alpha}})_{n\ge 1}$ is eventually log-concave for some $\alpha >3/2,$
 
 (ii) the sequence $\left ({\log (\sigma(-n))\over n}log(n)^A\right )_{n\ge 1}$ is eventually nondecreasing for some $A>0,$
 
 (iii) $\sum \limits_{n=1}^{+\infty}{log(\sigma(-n))\over n^2}=+\infty,$
 
 (iv) $\limsup_{n\to +\infty}{log \Vert T_{|_E}^n\Vert \over log(\sigma(-n))}<{1\over 200}.$
 
 Then for every closed left-invariant subspace $N\neq \{0\}$ of $F=E\oplus H^2_{0, \sigma}$  there exists $k\ge 0$ such that ${\bf S}^k(N)$ is analytic, and we have $F={\bf S}^k(N)+E$ and $N=\left [ \cup_{n\ge k}{\bf S}^{-n}\left ( {\bf S}^{k}(N)\cap E\right ) \right ]^-.$  In particular every nontrivial closed bi-invariant subspace $N$ of $F$ has the form $N=\left [ \cup_{n\ge k}{\bf S}^{-n}\left ( M\cap E\right ) \right ]^-,$ where $M=N\cap E$ is a closed $z$-invariant subspace of $E$ having the division property.
 
 \end{thm}
 
 Proof: Since $E$ and $\sigma$ satisfy the hypothesis of Proposition 7.5, and since the assertion concerning closed bi-invariant subspaces of $F$ is a consequence of the assertion concerning closed left-invariant subsaces, it suffices to show that for every closed left-invariant subspace $N\neq \{0\}$ of $F=E\oplus H^2_{0, \sigma}$  there exists $k\ge 0$ such that ${\bf S}^k(N)$ is analytic.

 It follows from (i), see Proposition 2.5 of \cite{ev1}, which is a discrete version of a result from \cite{bh}, that there exist a function $\phi \in W(\sigma)$ such that for every $\delta \in (0, {3\over 2 \alpha})$ there exists $k_{\delta}>0 $ satisfying, for every $h \in H^2_{0,\sigma}(\C \setminus \overline{\DD}).$
 
 \begin{equation}\left | \overline{\partial}\mathcal D_{\phi}(h) (\lambda)\right  | \le k_{\delta}\Vert h \Vert L^{-\delta}(\vert \lambda \vert) \ \ (\vert \lambda \vert <1),\end{equation}.
 
 where ${\mathcal D}_{\phi}(h)$ denotes the Dynkin extension of $w$ to $\DD$ associated to $\phi,$ and where $L(r)=\sup_{n\ge 0}{r^n\over \sigma^*(n)}=\sup_{n\ge 0}r^n\sigma(-n-1)$ is the so-called Legendre transform of the weight $\sigma^*.$

We will say that $l=(v, \theta)\in E^*\oplus H^2_{0, \sigma^*}(\C \setminus \overline{\DD})=F^*$ is not left-cyclic if there exists $f=(g,h)\in F\setminus\{0\}$ such that $\langle {\bf S}^{-n}f, l\rangle=0$ for every $n\ge 1.$ We now show that there exists $s_0\in (0, 1/4)$ such that $\sum \limits_{n=0}^{+\infty}\vert \widehat {\theta}(n)\vert ^2\sigma^{-2s_0}(-n)<+\infty$ for every $l=(v,\theta)\in F^*$ which is not left-cyclic.

To see this, assume that $l=(v, \theta)\in F^*$ is not left-cyclic,  and let $f=(g,h) \in \F\setminus \{0\}$ such that $\langle {\bf S}^{-n}f,l\rangle=0$ for every $n \ge 1.$ 
 
 Set $U=g +\mathcal D_{\phi}(h),$ set $ 
 V=\mathcal C^+\left(\theta \overline {\partial}\mathcal D_{\phi}(h)\right ) \in L^1(\DD),$ and set $W:=-\langle g_., v\rangle.$
 
 It follows from Proposition 5.4 that we have
 
 $$\theta U=V+W.$$
 
 Then $U(\lambda)=O\left( \Delta_E \left (\vert \lambda \vert\right )\right )$ and $H(\lambda)=O\left( \Delta_E \left (\vert \lambda \vert\right )\right )$ as $\vert \lambda \vert \to 1^-,$ and equation (27) gives an estimate of the rate of decrease of $\left \vert \overline{\partial}U(\lambda)\right \vert$ as $\lambda  \to 1^-.$ The existence of a number $s_0\in (0,1/4),$ which does not depend on $l,$ such that $\sum \limits_{n=0}^{+\infty}\vert \widehat {\theta}(n)\vert ^2\sigma^{-2s_0}(-n)<+\infty$ follows then directly from Lemma 5.1 (ii), Lemma 5.2, Lemma 5.5, Lemma 5.6 and Corollary 5.7 of \cite{ev1} which remain valid in this context. without any modification 
 
 We can now conclude the proof using the same duality argument as in the prrof of Theorem 5.8 of \cite{ev1}. Let $N$ be a nontrivial closed left-invariant subspace of F, let $f\in N\setminus \{0\}$ and let $k\ge 0,$ $w \in \ho$ and $g \in E$ satisfying the conditions of Theorem 8.3 with respect to $f$ and $s_0.$
 
 Let $l \in F^*$ such that $\langle u , l\rangle=0$ for every $u \in {\bf S}^k(N).$ Then $l$ is not left-cyclic, and so $l\in F_{s_0}^*$, where $F_{s_0}:=E\oplus H^2_{0, \sigma^{s_0}}.$ But it follows from conditions (1) and (2) of theorem 8.3 that $g$ belongs to the closure of ${\bf S}^k(N)$ in $F_s,$ and so $\langle g, l\rangle=0,$ $g \in {\bf S}^k(M),$ and ${\bf S}^k(M)$ is an analytic left-invariant subspace of $F,$ which concludes the proof of the theorem. $\square$
 
 \smallskip
 
 Of course, we see a posteriori as in \cite{ev1} that if the conditions of theorem 8.5 are satisfied, and if $N$ is a nontrivial closed left-invariant subspace of $F,$ then the map $f +{\bf S}^{k}(N)\cap E \to f +N$ defines an isomorphism from the quotient space $E/\left ( {\bf S}^k(N)\cap E\right )$ onto $F/N$, and so $\Vert {\bf S}^{-n}_{F/N}\Vert =O\left (\Vert T_E^{n+k}\Vert \right )$ as $n\to +\infty,$ where  ${\bf S}^{-n}_{F/N}$ denotes the operator induced  by ${\bf S}_{|_F}^{-n}$ on $F/N.$ This gives a much better estimate on $\widehat \theta(n)$ than the fact that $\sum \limits_{n=0}^{+\infty} \vert \widehat \theta(n)\vert^2\sigma^{2s_0}(-n-1)<+\infty $ if $l=(v, \theta) \in N^{\perp}.$ In particular if we apply Theorem 8.3 when the stronger hypothesis of Theorem 8.5 are satisfied, the function $g \in E$ associated to ${\bf S}^kf$ by Theorem 8.3 generates the same closed left-invariant subspace as ${\bf S}^kf.$

\bigskip

\bigskip

IMB, UMR 5251

Universit\'e de Bordeaux

351, cours de la Lib\'eration

33405 - Talence (France)

{\it esterle@math.u-bordeaux.fr}
\end{document}